\documentclass[12pt,a4paper]{scrartcl}
\usepackage{latexsym}
\usepackage{amsmath,amsfonts,amssymb,amsthm}
\usepackage{graphicx}
\usepackage{xcolor}
\usepackage{tikz}
\usepackage{booktabs}
\usepackage{tabularx}
\usepackage{xspace}
\usepackage{url}
\usepackage{paralist}

\makeatletter
\renewcommand\paragraph{\@startsection{paragraph}{4}{\z@}%
  {-3.25ex\@plus -1ex \@minus -.2ex}%
  {1.5ex \@plus .2ex}%
  {\normalfont\normalsize\bfseries}}
\makeatother

\usepackage{listings}

\lstdefinelanguage{Perl}
{morekeywords={rule,property,user_method,Int,Matrix,Vector,Set,Array,Integer, Rational,new,Polytope,row,rows,minor,col,cols,zero_vector,contains,declare, object,this,abs,accept,alarm,atan2,bind,binmode,bless,caller,%
      chdir,chmod,chomp,chop,chown,chr,chroot,close,closedir,connect,%
      continue,cos,crypt,dbmclose,dbmopen,defined,delete,die,do,dump,%
      each,else,elsif,endgrent,endhostent,endnetent,endprotoent,%
      endpwent,endservent,eof,eval,exec,exists,exit,exp,fcntl,fileno,%
      flock,for,foreach,fork,format,formline,getc,getgrent,getgrgid,%
      getgrnam,gethostbyaddr,gethostbyname,gethostent,getlogin,%
      getnetbyaddr,getnetbyname,getnetent,getpeername,getpgrp,%
      getppid,getpriority,getprotobyname,getprotobynumber,getprotoent,%
      getpwent,getpwnam,getpwuid,getservbyname,getservbyport,%
      getservent,getsockname,getsockopt,glob,gmtime,goto,grep,hex,if,%
      import,index,int,ioctl,join,keys,kill,last,lc,lcfirst,length,%
      link,listen,local,localtime,log,lstat,m,map,mkdir,msgctl,msgget,%
      msgrcv,msgsnd,my,next,no,oct,open,opendir,ord,pack,package,pipe,%
      pop,pos,print,printf,prototype,push,q,qq,quotemeta,qw,qx,rand,%
      read,readdir,readlink,recv,redo,ref,rename,require,reset,return,%
      reverse,rewinddir,rindex,rmdir,s,scalar,seek,seekdir,select,%
      semctl,semget,semop,send,setgrent,sethostent,setnetent,setpgrp,%
      setpriority,setprotoent,setpwent,setservent,setsockopt,shift,%
      shmctl,shmget,shmread,shmwrite,shutdown,sin,sleep,socket,%
      socketpair,sort,splice,split,sprintf,sqrt,srand,stat,study,sub,%
      substr,symlink,syscall,sysopen,sysread,system,syswrite,tell,%
      telldir,tie,tied,time,times,tr,truncate,uc,ucfirst,umask,undef,%
      unless,unlink,unpack,unshift,untie,until,use,utime,values,vec,%
      wait,waitpid,wantarray,warn,while,write,y},
morecomment=[l]{\#},
sensitive=true,
MoreSelectCharTable=\lst@ReplaceInput{\$\#}{\lst@ProcessOther\$\lst@ProcessOther\#},
}
\lstdefinelanguage{xml}
{morekeywords={property,value,name,object,type,m,v,?xml,version,encoding},
sensitive=true,
}

\providecommand\polymake{\texttt{polymake}\xspace}
\providecommand\normaliz{\texttt{normaliz2}\xspace}
\providecommand\lattE{\texttt{Latte macchiato}\xspace}
\providecommand\fourtitwo{\texttt{4ti2}\xspace}
\newcommand{\R}{\mathbb R}
\newcommand{\Q}{\mathbb Q}
\newcommand{\Z}{\mathbb Z}
\newcommand{\C}{\mathbb C}
\newcommand{\PP}{\mathbb P}

\newcommand{\pR}{\R_{\geq 0}}

\newcommand{\define}[2]{{\bf #1}\index{#2}}
\newcommand{\limit}{N}
\newcommand{\nlimit}{12 }
\newcommand{\scp}[1]{\langle #1 \rangle}

\newcommand{\edgeparam}[2]{
\fill[fill=black!0] (#1) circle (9pt);
\draw (#1) node[anchor=mid] {$#2$};
}
\newcommand{\edgeparamx}[2]{
\fill[fill=black!0] (#1) circle (10pt);
\draw (#1) node[anchor=mid,rotate=30] {\tiny$#2$};
}

\DeclareMathOperator{\inter}{int}
\DeclareMathOperator{\cone}{cone}
\DeclareMathOperator{\conv}{conv}
\DeclareMathOperator{\spann}{span}
\DeclareMathOperator{\RHS}{RHS}
\DeclareMathOperator{\NF}{NF}
\DeclareMathOperator{\NC}{NC}
\DeclareMathOperator{\nb}{nb}
\DeclareMathOperator{\B}{B}

\newtheorem*{theorem}{Theorem}
\newtheorem*{prop}{Proposition}

\theoremstyle{definition}

\newtheorem*{defn}{Definition}

\theoremstyle{remark}
\newtheorem*{lem}{Lemma}
\newtheorem*{pf}{Proof}
\newtheorem*{rem}{Remark}

\title{Classification of smooth lattice polytopes with few lattice points}
\author{Benjamin Lorenz}

\begin{document}

\begin{titlepage}
\begin{center}
\Huge 
{\bf Freie Universit\"at Berlin}\\
{\bf Fachbereich Mathematik}\\\vspace{2cm}
{\bf Diplomarbeit}\\\vspace{2cm}
\LARGE Classification of smooth lattice polytopes\\ with few lattice points\\\vspace{1.5cm}
\large Supervisor: Dr. Christian Haase
\end{center}\vspace{1cm}
{
{\bfseries Abstract:}
After giving a short introduction on smooth lattice polytopes, I will present a proof for the finiteness of smooth lattice polytopes with few lattice points. The argument is then turned into an algorithm for the classification of smooth lattice polytopes in fixed dimension with an upper bound on the number of lattice points. Additionally I have implemented this algorithm for dimension two and three and used it, together with a classification of smooth minimal fans by Tadao Oda, to create lists of all smooth 2\nobreakdash-polytopes and 3\nobreakdash-polytopes with at most 12 lattice points.
}\\
\vspace{5cm}\\
\large
\begin{tabular}{ll}
by: & Benjamin Lorenz\\
& lorenz.benjamin@fu-berlin.de
\end{tabular}
\normalsize
\thispagestyle{empty}
\end{titlepage}

\tableofcontents

\pagebreak[4]

\section{Introduction}
\subsection{Polytopes and Lattices}
\par
Our main goal is the classification of smooth lattice polytopes with few lattice points, so we will start by introducing polytopes and lattice polytopes. 

\begin{defn}
A set $P\subset\R^d$ is called a \define{polytope}{polytope} if it is the convex hull of a finite set of points $V = \{v_1,\ldots,v_n\}\subset\R^d$. If $V$ is minimal up to inclusion, then the points $v_1,\ldots,v_n$ are the vertices of $P$.
\[
P = \conv(V) = \left\{ \sum_{v\in V}\lambda_v v\in\R^d \ |\ \sum_{v\in V} \lambda_v = 1, \lambda_v\ge 0\ \forall v\in V \right\}
\]
\end{defn}
\begin{rem}
In the following we will always assume $V$ to be the vertex set of the polytope.
\end{rem}
\par
For the definition of smooth lattice polytopes, and to make precise what we mean by polytopes with few lattice points, we first need to define lattices in $\R^d$.
\begin{defn}
For some finite set $B=\{b_1,\ldots,b_n\}\subset \R^d$, we define the \define{lattice}{lattice} $\Lambda(B)$: 
\[
\Lambda(B) = \left\{\sum_{b\in B} \lambda_b b\in\R^d\ |\ \lambda_b \in \Z \ \forall b\in B\right\}\subset\R^d
\]
If the vectors in $B$ are linearly independent then $B$ is called a \define{basis}{lattice!basis} of $\Lambda(B)$.
\end{defn}
\begin{rem}
We will now restrict to the lattice $\Lambda = \Z^d\subset\R^d$, so any set $B=\{b_1,\ldots,b_d\}\subset \Z^d$ is a basis for $\Z^d$ if $\det(B) = \pm 1$.
\end{rem}
Now, for a lattice $\Lambda=\Lambda(B)\subset\R^d$, the \define{dual lattice}{dual!lattice} is
\[
\Lambda^*=\{\gamma\in(\R^d)^*|\ \gamma(x)\in\Z \ \forall x\in\Lambda\}
\]
and one can easily see that for $\Lambda(B)=\Z^d$ we have $\Lambda^*=\Lambda((B^{-1})^T)\cong\Z^d$. So for some lattice basis $B$ define the \define{dual lattice basis}{lattice!dual basis} $B^*=(B^{-1})^T$.

\par
We can now define smooth lattice polytopes:
\begin{defn}
Let $P$ be a polytope. Then $P$ is called a \define{lattice}{polytope!lattice} polytope if all vertices have integer coordinates, i.e.~$V \subseteq \Z^d$.
\end{defn}
When we talk about the lattice points $L(P)$ of a polytope $P\subset\R^d$ we mean the elements of $\Z^d$ that lie in the polytope, namely $L(P) = P\cap\Z^d$. So we want to limit the number $|L(P)|$ for our polytopes.
\par
For a lattice polytope $P$ to be smooth it has to be simple. So for $v\in V$ denote by $\nb(v)$ the set of adjacent vertices, i.e.~for every $w\in\nb(v)$ there is an edge containing $v$ and $w$. Then a polytope $P\subset\R^d$ is \define{simple}{polytope!simple} if $|\nb(v)| = d$ for all vertices $v\in V$.
\begin{defn}
A simple lattice polytope $P$ is called \define{smooth}{polytope!smooth} if, at each vertex, the lattice minimal edge-directions form a lattice basis, or more formally
\[
\det\left(\left\{\frac{w-v}{\gcd(w-v)} |\ w\in \nb(v) \right\}\right) = \pm 1 \quad \forall v \in V
\]
\end{defn}

\subsection{Cones and Fans}
In our main algorithm we want to approach the polytopes from the dual and therefore we now introduce cones and fans and define some important properties.
\begin{defn}
For $S\subset \R^d$ finite, the \define{polyhedral cone}{cone} $C=\cone(S)$ is 
\[
\cone(S) = \left\{\sum_{r\in S} a_r r\in\R^d |\ a_r \in \pR \ \forall r\in S \right\}
\]
A cone $C\subset\R^d$ is called \define{pointed}{cone!pointed} if there is $\alpha\in(\R^d)^*$ such that 
\[
C\cap\{x\in\R^d\ |\ \alpha(x) \le 0\} = \{0\}
\]
If $C=\cone(S)$ is pointed and $S$ is minimal up to inclusion, then the vectors $r \in S$ are called the \define{rays}{cone!ray} of $C$.
\end{defn}
\begin{rem}
Some restrictions on the cones we will consider later on.
\begin{itemize}
\item All our cones are pointed.
\item All our cones are rational, which means that there is $S\subset\Q^d$ finite, such that $C=\cone(S)$.
\item And we assume all rays $r\in S$ to be primitive vectors, i.e.~$\gcd(r)=1$ for all $r\in S$.
\end{itemize}
\end{rem}
The \define{dimension}{cone!dimension} of a cone $C=\cone(S)$ is the dimension of its linear hull $\dim(C)=\dim(\spann(S))$. A $k$-dimensional cone $C$ is called \define{simplicial}{cone!simplicial} if there is $S=\{r_1,\ldots,r_k\}$ such that $C=\cone(S)$. One may already see the corresondence to simple polytopes. We will continue in this direction by defining smooth cones. Later we will see that both definitions of smoothness coincide if we choose the right construction to define a fan from a polytope.
\begin{defn}
Let $C\subset\R^d$ be a $k$-dimensional simplicial cone, then $C$ is \define{smooth}{cone!smooth} if there are $z_1,\ldots,z_d\in\Z^d$ such that
\begin{itemize}
\item $\ C=\cone(\{z_1,\ldots,z_k\}) $
\item $\ \det(z_1,\ldots,z_d)=\pm 1 $
\end{itemize}
\end{defn}
\pagebreak[3]
\begin{defn}
A \define{fan}{fan} $F$ is a finite set of cones with the following properties.
\begin{enumerate}
\item[(1)] $C_1,C_2\in F \Rightarrow C_1\cap C_2 \in F$
\item[(2)] $C \in F,\ C'\ne\emptyset,\ C' \le C$ (i.e.~$C'$ is a non-empty face of $C$) $ \Rightarrow C'\in F$
\end{enumerate}
\end{defn}
For a fan $F$, denote by $F^{(k)}$ the set of cones of dimension $k$.
A fan $F$ in $\R^d$ is called \define{complete}{fan!complete} if it covers the whole $\R^d$, which means $\bigcup_{C\in F}C = \R^d$.\\
The properties of cones defined above can of course be extended to fans, a fan  is \define{simplicial}{fan!simplicial} or \define{smooth}{fan!smooth} if all its cones are.

\subsection{Duality}
In the third section of the introduction we will now state a duality theorem of polyhedral theory and establish a connection between polytopes and fans. This will turn out useful later in the classification.
\begin{defn}
For a matrix $A\in\R^{m\times d}$ and a vector $b\in\R^m$ the \define{polyhedron}{polyhedron} $P(A,b)\subseteq\R^d$ is the intersection of finitely many halfspaces
\[
P(A,b) = \bigcap_{a_i \mbox{\tiny\ row of } A} \{x\in\R^d |\ \scp{a_i,x}\le b_i\} = \{ x\in\R^d |\ Ax \le b\}
\]
\end{defn}
For the duality theorem we need the \define{Minkowski sum}{Minkowski sum} $P+Q$ of two subsets $P,Q\subset\R^d$:
\[
P+Q= \{ x+y\in\R^d \ | \ x\in P, \ y\in Q\}
\]
\begin{theorem}
A set $P\subseteq\R^d$ is a polyhedron if and only if it is the Minkowski sum of a polytope and a cone. 
\end{theorem}
\begin{pf}
The proof of this duality can be found in most textbooks containing polyhedral theory, e.g. {\itshape Lectures on Polytopes}~\cite{ziegler}.
\end{pf}
\begin{rem}
In particular we have the following results for polytopes:
\begin{itemize}
\item A polytope is a bounded polyhedron.
\item For $P=P(A,b)$ bounded, there is $V\subset\R^d$ finite, such that $P=\conv(V)$.
\item And for $P=\conv(V)$ there are $A\in\R^{m\times d}$ and $b\in\R^m$ such that $P=P(A,b)$.
\end{itemize}
\end{rem} 
This is very useful as it allows us to switch between both representations as we need it. Now let us define the normal fan of a polytope, the main construction we will use for the classification of the polytopes.
\begin{defn}
Let $P=P(A,b)\subset\R^d$ be a polytope, the \define{normal cone}{normal cone} of a face $F\subset P$ is
\[
\NC_P(F)=\cone(\{a_i\mbox{ row of }A\ | \ \scp{a_i,x}=b_i \ \forall x\in F\})
\]
The \define{normal fan}{normal fan} of $P$ is defined by the normal cones of all faces of the polytope.
\[
\NF_P = \{\NC_P(F) \ | \ F \subset P \mbox{ face} \} 
\] 
\end{defn}

To see that this does indeed produce a fan, observe that for every face $F$ of the polytope the faces containing $F$ are defined by subsets of the rays defining $F$. Thus, the faces of $P$ containing $F$ define the faces of $\NC_P(F)$. Additionally this construction always produces a complete fan.\\
\begin{minipage}[b]{0.8\textwidth}
\vspace{2mm} This yields a construction to build a fan from a given polytope. We can also use this in the other direction. For a given matrix $A$ containing the rays of the fan we can define a polytope $P(A,b)$ by choosing a right-hand side vector $b$. We just need to make sure that the normal fan of the polytope is again the fan we started with, since there can be different polytopes with the same set of facet normals.
\end{minipage}
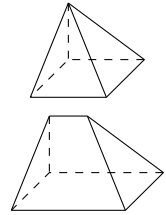
\begin{minipage}[b]{0.2\textwidth}
\begin{center} \begin{tikzpicture}[scale=0.5]
\draw (0,0) -- (2,0) -- (3,1);
\draw[dashed] (3,1) -- (1,1) -- (0,0);
\draw (1,2.5) -- (0,0);
\draw (1,2.5) -- (2,0);
\draw (1,2.5) -- (3,1);
\draw[dashed] (1,2.5) -- (1,1);
\end{tikzpicture}\\\vspace{2mm}
\begin{tikzpicture}[scale=0.5]
\draw (0,0) -- (3,0) -- (4,1);
\draw[dashed] (4,1) -- (1,1) -- (0,0);
\draw (1,2.5) -- (0,0);
\draw (2,2.5) -- (3,0);
\draw (2,2.5) -- (4,1);
\draw[dashed] (1,2.5) -- (1,1);
\draw (1,2.5) -- (2,2.5);
\end{tikzpicture}
\end{center}
\end{minipage}\\
However, we only want to classify smooth polytopes, so the obvious question is whether there is a property of the normal fan that coincides with the smoothness of the polytope. And the obvious answer is also correct:
\begin{prop}
Let $P$ be lattice $d$-polytope and $\NF_P$ its normal fan. Then $P$ is smooth if and only if $\NF_P$ is smooth.
\end{prop}
\begin{pf}
First we will show that $P$ is simple if and only if $\NF_P$ is simplicial.
\begin{itemize}
\item Let $\NF_P$ be a simplicial fan, then the normal cone of a vertex $v$ of $P$ is defined by exactly $d$ rays. Each $(d-1)$ subset of the facets defined by these rays can define an edge and since there are $d$ such subsets, there can be at most $d$ edges containing $v$. Because $P$ is a $d$-dimensional polytope, there must be at least $d$ edges containing $v$ and thus $P$ is simple.
\item Let $P$ be a simple polytope and $v$ be a vertex of $P$. Now, every facet containing $v$ is itself a $(d-1)$-dimensional polytope and thus there must be at least $(d-1)$ edges containing $v$ in each facet. Since there are only $d$ $(d-1)$-subsets of the $d$ edges at $v$ there can be at most $d$ different facets and thus $\NF_P$ is simplicial.
\end{itemize}
We choose again one vertex $v$ in $P$. To prove the smoothness we use the matrizes $E,\ A$ defined by the edge-directions $e_1,\ldots,e_d$ at $v$ and the normal cone of the vertex $N_P(v) = \cone(\{a_1,\ldots,a_d\})$.\\
For given smooth edge-directions $E$ choose $A=-(E^{-1})^T$, obviously $A$ defines the normal cone. By Cramer's rule and since $\det(E)=\pm 1$ we know that $A$ is integral and $\det(A)=\mp 1$, thus the cone is smooth. And for given smooth normal fan with rays in $A$ choose $E=-(A^{-1})^T$. By the same argument $E\in\Z^{d\times d}$, $\det(E)=\mp 1$ and thus $E$ contains the lattice minimal edge-directions.
\qed
\end{pf}

\pagebreak
\section{Finiteness}
\subsection{Finiteness Theorem}
Now we can state the main theorem~\cite{finiteness}.
\begin{theorem}
For integers $d$ and $\limit$, there are only finitely many smooth
$d$-polytopes with at most $\limit$ lattice points.
\end{theorem}
Since the following proof inspired the technique for the classification algorithm following subsequently, we will repeat its details.
\begin{pf}
We will proof this in three steps. First we want to limit the number of combinatorial types of smooth fans, then the number of smooth fans with each of those types and finally the number of smooth polytopes having one of those as normal fan, all using the limit on the number of lattice points.\\
In the following we will again assume all rays to be primitive.
\begin{asparaitem}
\item Assume we have a combinatorial type of a complete simplicial fan $F$ with $n$ rays and at most $N$ full-dimensional cones. A full-dimensional cone of this fan can be defined by $d$-subsets $C \subset \{1,\ldots,n\},\ |C|=d$. The fan can then be defined by the set of all such full-dimensional cones $F^{(d)}\subset\binom{\{1,\ldots,n\}}{d}$. {\itshape Note that we use a different notation here since we consider only the combinatorial type of the fan and do not care about the rays.} We know that the fan has at most $N$ full-dimensional cones, so $|F^{(d)}| \le N$. Additionally every ray needs to appear in at least $d$ full-dimensional cones:
\[
d_i= |\{C\in F^{(d)} \ | \ i\in C\}| \ge d \quad \quad \forall i\in\{1,\ldots,n\}
\]
By counting all appearances of all rays in $F^{(d)}$ in two different ways we have:
\[
d\cdot n \le \sum_{i=1}^n d_i \le d\cdot N
\]
And thus $n\le N$, which means that we have at most $N$ rays in the simplicial fan $F$. Hence our fan is fully defined by choosing 
\[
F^{(d)}\subset\binom{\{1,\ldots,n\}}{d}\subset\binom{\{1,\ldots,N\}}{d}
\]
Since $N$ and $d$ are fixed there are only finitely many such subsets, and thus only finitely many possible combinatorial types of simplicial fans.
\item Given the combinatorial type we now want to define the fan $F$ up to isomorphism:\\
Let $C=\cone(n_1,n_2,\ldots,n_{d-1})\in F$ be a codimension $1$ cone which is contained in two full-dimensional cones $C_1=\cone(C,r_1)$ and $C_2=\cone(C,r_2)$. Since all the cones are smooth $r_1$ and $r_2$ have lattice distance one from $C$, but lie on different sides of $C$. So $r_1+r_2$ lies in the span of $C$ and thus
\[
r_1 + r_2 = \sum_{i=1}^{d-1} a_{C,i} n_i
\]
The coefficients $a_{C,i}$ are the edge-paremeters of a fan $F$. We know $a_{C,i}\in\Z$ because $C$ is smooth and thus $\{n_1,\ldots,n_{d-1}\}$ is a Hilbert basis of $C$.\\
We can also use the edge-parameters $a_{C,i}$ for $1\le i\le d-1,\ C\in F^{(d-1)}$ together with the combinatorial type to reconstruct fan $F$ using the following technique: Every smooth fan can, while preserving the lattice, be transformed to contain the cone defined by all $d$ unit vectors. This means we can assume the fan to contain this cone. Then there are codimension $1$ cones for every $(d-1)$-subset of the unit vectors, so, given the edge-parameters for each one of those cones, one can calculate the ray on the other side of this codimension $1$ cone with the formula above. Iterate until all rays are calculated.\\
Since the edge-parameters are in $\Z$ there may be infinitely many such smooth fans, but in the next step we will show a way to bound those parameters.

\item In the last step we have to bound the number of smooth polytopes $P$ with at most $N$ lattice points having a given normal fan $F$ and the number of smooth fans the above construction yields.\\ 
We therefore consider the edge $e$ of the polytope dual to a codimension $1$ cone $C\in F$. Let $l(e)$ be the lattice length of the edge $e$ and $x_1, x_2$ be the vertices defining this edge. Denote by $u_1,\ldots,u_d$ the lattice-minimal edge-directions at $x_1$ which are dual to the inner normals $-n_1,\ldots,-n_{d-1},-r_1$ around $x_1$. So we have $x_2 = x_1 + l(e) u_d$. Using the edge-directions at $x_1$ and the edge-parameters we can calculate the edge-directions at $x_2$.
\begin{lem}
For $x_1,x_2$ as above, the edge-directions at $x_2$ are
\begin{align*}
u_i'&=u_i+a_{C,i}u_d \quad \mbox{ for } i\in\{1,\ldots,d-1\}\\
u_d'&=-u_d
\end{align*}
\end{lem}
\begin{pf}
The inner normals around $x_2$ are $-n_1,\ldots,-n_{d-1},-r_2$ and 
\[
-r_2 = r_1 - \sum_{i=1}^{d-1} a_{C,i} n_i
\] so we have
\begin{align*}
\scp{-n_i,u_j'} & =\scp{-n_i,u_j+a_{C,j}u_d} =\scp{-n_i,u_j}-a_{C,j}\scp{n_i,u_d} & \mbox{ for } i, j < d \\
 & = \scp{-n_i,u_j} = \delta_{i,j} & \mbox{ for } i, j < d \\
\scp{-n_i,u_d'} & = \scp{-n_i,-u_d} = 0 & \mbox{ for } i < d
\end{align*}
\begin{align*}
\scp{-r_2,u_j'} & = \scp{r_1,u_j+a_{C,j}u_d} - \sum_{i=1}^{d-1} a_{C,i} \scp{n_i,u_j+a_{C,j}u_d} & \mbox{ for } j < d \\
 & = a_{C,j}\scp{r_1,u_d} - a_{C,j}\scp{n_j,u_j} = 0 & \mbox{ for } j < d \\
\scp{-r_2,u_d'} & = \scp{r_1,-u_d} - \sum_{i=1}^{d-1} a_{C,i} \scp{n_i,-u_d} = \scp{-r_1,u_d} = 1 & 
\end{align*}
which proves that the matrix containing these edge-directions is inverse to the given inner normals, finishing the proof of the lemma.
\qed
\end{pf}
Since in our smooth lattice polytope every edge must have length at least $1$ the lattice points defined by the vertex $x_2$ plus one of the edge-directions must of course lie in the polytope and thus 
\begin{align*}
\scp{r_1,x_2+a_{C,i}u_d} & \ge \scp{r_1,x_1} \\
\scp{r_1,x_1+l(e)u_d+a_{C,i}u_d} & \ge \scp{r_1,x_1} \\
(l(e)+a_{C,i})\scp{r_1,u_d} &\ge 0 \\
a_{C,i} &\ge -l(e) \ge -N
\end{align*}
This gives a lower bound for the edge-parameters. For an upper bound let us first introduce the thickened edge.
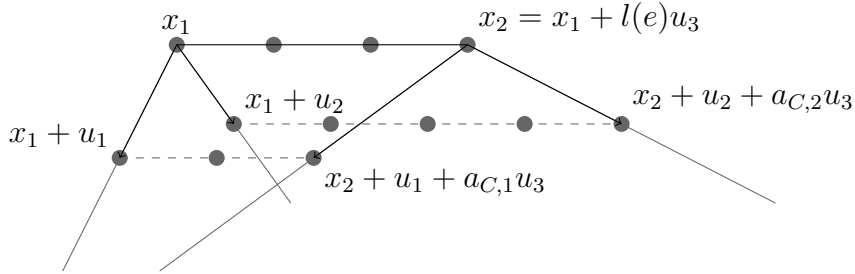
\begin{figure}[htb]
\begin{center}
\begin{tikzpicture}[scale=1.5]
\draw[black!60] (0,0) -- (1,2);
\draw[black!60] (2,0.6) -- (1,2);
\draw[black!60,dashed] (0.5,1) -- (2.2,1);
\draw[black!60,dashed] (1.5,1.3) -- (4.9,1.3);
\draw[black!60] (3.55,2) -- (0.85,0);
\draw[black!60] (3.55,2) -- (6.25,0.6);

\fill[black!60] (0.5,1) circle (2pt);
\fill[black!60] (1.5,1.3) circle (2pt);
\fill[black!60] (1,2) circle (2pt);
\fill[black!60] (1.85,2) circle (2pt);
\fill[black!60] (2.7,2) circle (2pt);
\fill[black!60] (3.55,2) circle (2pt);
\fill[black!60] (1.35,1) circle (2pt);
\fill[black!60] (2.2,1) circle (2pt);
\fill[black!60] (2.35,1.3) circle (2pt);
\fill[black!60] (3.2,1.3) circle (2pt);
\fill[black!60] (4.05,1.3) circle (2pt);
\fill[black!60] (4.9,1.3) circle (2pt);

\draw[black,->] (3.55,2) -- (2.2,1) node[anchor=north west] {$x_2+u_1+a_{C,1}u_3$};
\draw[black,->] (3.55,2) -- (4.9,1.3) node[anchor=south west] {$x_2+u_2+a_{C,2}u_3$};
\draw[black,->] (1,2) -- (0.5,1) node[anchor=south east] {$x_1+u_1$};
\draw[black,->] (1,2) -- (1.5,1.3) node[anchor=south west] {$x_1+u_2$};
\draw[black] (1,2) node[anchor=south] {$x_1$} -- (3.55,2) node[anchor=south west] {$x_2=x_1+l(e)u_3$};
\end{tikzpicture}
\end{center}
\caption{A thickened edge}
\label{thickenededge}
\end{figure}
\begin{defn}
Using the notation from above we define the thickened edge in a smooth polytope as the convex hull of all those points we just discussed.
\begin{align*}
c(e) = \conv(\{&x_1,x_1+u_1,\ldots,x_1+u_{d-1}, \\
&x_2,x_2+u_1+a_{C,1}u_d,\ldots,x_2+u_{d-1}+a_{C,d-1}u_d \})
\end{align*}
\end{defn}
We obviously have $c(e)\subseteq P$ and as shown in Figure~\ref{thickenededge} $c(e)$ contains $d-1$ lines from $x_1+u_i$ to $x_2+u_i+a_{C,i}u_d$ each containing $l(e)+1 + a_{C,i}$ lattice points and the line from $x_1$ to $x_2$ containing $l(e)+1$ lattice points, so 
\[
|c(e)\cap\Z^d| = d(l(e)+1) + \sum_{i=1}^{d-1}a_{C,i}
\]
Since the edge length must be positive and we want the number of lattice points of $P$ and thus of $c(e)$ to be less than $N$, this gives us an upper bound for the sum of all edge-parameters. Together with the lower bound this leaves us with only finitely many choices for the edge-parameters. This means for an upper bound $N$ on the number of lattice points in the polytope there are only finitely many possible smooth normal fans.\\
Now for a given normal fan we can define the polytope up to isomorphism by fixing all edge lengths. Just assume one vertex at the origin and the edge-directions to be the unit vectors, then the edge lengths define the next vertices and the normal fan defines the next edge-directions. Iterate until all vertices are fixed. And we want of course the edge length $l(e)$ to be positive and less than $N$, so for every edge there are only finitely many possible edge lengths. And thus only finitely many polytopes for this fan.
\end{asparaitem}
So altogether we have finitely many combinatorial types of fans, finitely many choices for the parameters and finitely many choices for the edge lengths and thus only finitely many smooth lattice $d$-polytopes with at most $N$ lattice points. This completes the proof of the finiteness theorem.
\qed
\end{pf}

\subsection{General Algorithm}
The proof gives rise to the following approach for the classification of such polytopes.

\subsubsection{Generating Normal Fans}
The first step is to find all combinatorial types of complete simplicial fans in dimension $d$ with at most $N$ full-dimensional cones. One can then proceed by finding parameters $a_{C,i}$ which produce smooth fans, for all codimension $1$ cones.\\
In dimension three we can use that the combinatorial types of simplicial fans are equivalent to the combinatorial types of triangulations of $S^2$ with at most $\limit$ $2$-simplices, which have been classified by R. Bowen and S. Fisk in \emph{Generation of triangulations of the sphere}~\cite{sphere}.
And Tadao Oda has done the classification of smooth two-dimensional fans which are minimal up to equivariant blowing-ups and the same in dimension three for at most eight rays in \emph{Convex Bodies and Algebraic Geometry}~\cite{oda}. We will explain and use those results later in chapter $3$.

\subsubsection{Generating Polytopes}
The next step now is to find the polytopes. In the proof we used the edge lengths together with the normal fan to define the polytope. Here we will use the right-hand side of its inequality description $P(A,b) = \{x\ | \ Ax\le b\}$ instead of the edge lengths, but after expressing the edge length in terms of the $b$-vector we can use the edge lengths to find all those $b$-vectors defining smooth lattice polytopes with the correct normal fan.

\paragraph{Right-Hand Side}
For a given fan $F$ we now need to find all possible $b$-vectors which yield smooth lattice polytopes with at most $\limit$ lattice points and $F$ as normal fan. Therefore define the set of all such vectors 
\[
\RHS(F,N) = \{ b \in \Z^m |\ \NF_{P(A,b)} = F,\ |P(A,b)\cap\Z^d|\le \limit \} 
\]
where the rows of $A$ are the rays of $F$. In the following we want to define a polytope containing this set and enumerate all its lattice points.

\paragraph{Calculating the edge lengths}
\begin{lem}
Let $C=\cone(n_1,n_2,\ldots,n_{d-1})$ be a codimension $1$ cone and $r_1, r_2\in\Z^d$ such that $ C_1=\cone(C,r_1),\ C_2=\cone(C,r_2) \in F $ are the full-dimensional cones containing $C$. Then the length $l(e)$ of the edge $e$ in the polytope $P(A,b)$ dual to $C$ can be expressed in terms of the right-hand side $b$ as follows:
\[
l(e) = b_{r_1} + b_{r_2} - \sum_{i=1}^{d-1} a_{C,i} b_{n_i}
\]
where $a_{C,i}\ (1\le i \le d-1)$ such that $r_1+r_2 = \sum_{i=1}^{d-1} a_{C,i} n_i$.
\end{lem}
\begin{pf}
The coefficients $a_*$ are the same as in the proof of the theorem. Now let $x_1$ and $x_2$ be the vertices dual to $C_1$ and $C_2$ respectively. From the proof of the finiteness theorem we know $\scp{u_d,-r_1}=1$ and thus $l(e) = \scp{x_2-x_1,-r_1}$. Since $r_1+r_2 = \sum_{i=1}^{d-1} a_{C,i} n_i$, this yields
\begin{align*}
l(e) & = \scp{x_1,r_1} - \scp{x_2,r_1}  \\
  & = b_{r_1} - \scp{x_2,\sum_{i=1}^{d-1} a_{C,i} n_i - r_2}  \\
  & = b_{r_1} - \sum_{i=1}^{d-1} a_{C,i} \scp{x_2,n_i} + \scp{x_2,r_2} \\
  & = b_{r_1} - \sum_{i=1}^{d-1} a_{C,i} b_{n_i} + b_{r_2}
\end{align*}
\qed
\end{pf}

\paragraph{Approximating the Right-Hand Sides}
Since we can express $l(e)$ linearly in terms of the $b$-vector we can now use bounds on the edge length to approximate $\RHS(F,N)$. For a codimension $1$ cone $C$, denote by $e_C$ the edge in the polytope dual to $C$.
\begin{lem}
For a given normal fan $F$ in dimension $d$ and a maximal number of lattice points $N$, every polytope defined by $F$ and a vector from $\RHS(F,N)$ satisfies the following inequalities
\begin{align*}
l(e_C) &\le \frac{\limit-\sum_{i=1}^{d-1} a_{C,i}}{d} - 1  & \forall C\in F^{(d-1)}\\
l(e_C) &\ge 1 & \forall C\in F^{(d-1)}\\
b_i & = 0 & \forall 1 \le i \le d \\
\sum_{C\in F^{(d-1)}}(l(e_C)-1) &\le N - |F^{(d)}| &
\end{align*}
Using $l(e)$ as defined in the previous lemma, these inequalities define a polytope \[\B(F,N)\supset\RHS(F,N)\]
\end{lem}
\begin{pf} We first show that $\RHS(F,N)$ satisfies the inequalities:
\begin{enumerate}
\item[(i)] We can now use an idea from the proof of the finiteness theorem. The convex hull of the thickened edge $c(e)$ contains $d(l(e)+1) + \sum_{i=1}^{d-1} a_{C,i}$ lattice points for an edge $e$ with the parameters $a_{C,i}$, so we have the following bound for $l(e)$.
\[
l(e) \le \frac{\limit-\sum_{i=1}^{d-1} a_{C,i}}{d} - 1
\]
\item[(ii)] The fan completely defines the combinatorial type of the polytope, so every codimension $1$ cone must define a non-degenerate (i.e.~$l(e) > 0$) edge in the polytope and since all vertices have integer coordinates the edge length must be at least one. This gives the second set of inequalities.
\item[(iii)] We now assume $F$ to be lattice-transformed such that the first $d$ rays are the unit vectors, this is possible because $F$ is smooth. Since we do not care about translations, we can fix the first $d$ right-hand sides to be zero.
\item[(iv)] The total number of lattice points on all edges of the polytope is the number of relative interior lattice points of each edge ($l(e)-1$) plus the number of vertices, so this sum must be less than the given bound $\limit$ on the number of lattice points.

\end{enumerate}
So all inequalities hold for any $b \in \RHS(F,N)$, they are all linear in terms of $b$ and thus $\B(F,N)$ is a polyhedron.\\
Now to prove that this polyhedron is bounded we use a similar approach as we used to construct the fan from the edge-parameters. For every codimension $1$ cone we have the following inequalities
\[
1 \le l(e) = b_{r_1} + b_{r_2} - \sum_{i=1}^{d-1} a_{C,i} b_{n_i} \le \limit
\]
Starting with the full-dimensional cone defined by the unit-vectors, where all $b_i$ ($1\le i \le d$) are zero, we can always find codimension $1$ cones where only one of the $b_*$ in the inequalities above is not yet bounded and thus bound it. We can iterate this until all $b_*$ are bounded and thus $\B(F,N)$ is a polytope.
\qed
\end{pf}

\paragraph{Enumeration}
Now, given an inequality description of $\B(F,N)$ we can calculate all lattice vectors $b\in \B(F,N)\cap\Z^m$. Together with the matrix $A$, each one defines a smooth lattice polytope $P(A,b)$ with possibly too many lattice points. Thus, after removing all polytopes with more than $N$ lattice points, we obtain a list of smooth lattice polytopes with normal fan $F$ and at most $N$ lattice points. The last step is to remove all polytopes that are lattice-isomorphic to another one. All these computations can be done with \polymake~\cite{polymake} using interfaces to \fourtitwo~\cite{4ti2}, \lattE~\cite{latte} and \normaliz~\cite{normaliz}.

\pagebreak
\section{Classification}
After a short definition of blowing-ups, the sections 3.2 and 3.3 contain a detailed description of how to classify smooth two- and three-dimensional lattice polytopes with at most $\nlimit$ lattice points based on the technique described above.

\subsection{Blowing-ups}
We want to use the results of Tadao Oda who classified smooth fans in dimension two and three which are minimal up to equivariant blowing-ups. This means that all normal fans we need can be generated by successive equivariant blowing-ups of the fans given in {\itshape Convex Bodies and Algebraic Geometry} by Tadao Oda~\cite{oda}. But first we need to clarify what such a blowing-up is. In short, it means subdividing a cone and all cones containing it.
\begin{defn}
Let $\hat C$ be a cone in the fan $F$ and $r \in \hat C$. The blowing-up $b(C,r)$ of a cone $C$ at $r$ is 
\begin{align*}
b(C,r) & = \{ \cone(\{r,C'\})\ |\ C'\in F, C'\le C, r\notin C' \} & \mbox{if } r\in C\\
b(C,r) & = \{C\} & \mbox{if } r\notin C
\end{align*}
\parbox[t]{0.6\linewidth}{
where $C'\le C$ means that $C'$ is a face of $C$.\\
Now the \define{blowing-up}{fan!blowing-up} $b(F,r)$ of $F$ at $r$ is 
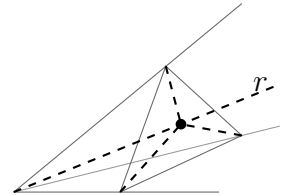
\[
b(F,r) = \bigcup_{C\in F}b(C,r)
\]
One can easily see that this is indeed again a fan.
}
\parbox[t]{0.4\linewidth}{
\[
\begin{tikzpicture}[scale=1.0]
\draw[black!60] (1,0) -- (3.7,0);
\draw[black!40] (1,0) -- (4.5,0.875);
\draw[black!60] (1,0) -- (4,2.5);
\draw[black!60] (2.4,0) -- (3,1.67) -- (4,0.75) -- cycle;
\draw[black,dashed,thick] (3,1.67) -- (3.2,0.9);
\draw[black,dashed,thick] (2.4,0) -- (3.2,0.9);
\draw[black,dashed,thick] (4,0.75) -- (3.2,0.9);
\draw[black,dashed,thick] (1,0) -- (4.5,1.43) node[anchor=east] {$r$};
\fill[black] (3.2,0.9) circle (2pt);
\end{tikzpicture}
\]
}
\end{defn}
A blowing-up is \define{equivariant}{blowing-up!equivariant} if $r=\sum_{s\in S}s$ for some cone  $C=\cone(S)\in F$, where $S$ is minimal up to inclusion and all rays are primitive. We denote by $b(F,C)$ the blowing-up in cone $C$.\\
For a smooth fan, an equivariant blowing-up is again smooth. Since we want to classify smooth polytopes and thus smooth normal fans, we only have to consider equivariant blowing-ups. \\
\begin{rem}
Since we have a connection between the fan and the polytope one could ask what such a blowing-up of the normal fan corresponds to on the polytope side. \\
We first choose a cone in the fan, which is equivalent to choosing a face in the polytope, and then add the sum of all rays of the cone as new ray. This corresponds to adding a new facet by cutting-off that face of the polytope.
\end{rem}

\subsection{Smooth Polytopes in Dimension Two}
\subsubsection{Classifying the Normal Fans}

As already mentioned we do not start with combinatorial types of fans, but with minimal smooth fans, and in dimension two there are only two such fans as given by the following theorem.

\paragraph{Minimal Fans}

\begin{theorem}{\bf(Tadao Oda~\cite{oda})}
Every smooth fan in dimension two can, up to isomorphism, be created by finitely many successive equivariant blowing-ups of one of the following fans.
\begin{itemize}
\item[(i)] The fan $F_p$ representing the complex projective space $\PP_2(\C)$.
\begin{align*}
F_p^{(1)}&= \{a_0,a_1,a_2\} = \{(1,0),(0,1),(-1,-1)\} \\
F_p^{(d)}&=\{\cone(\{a_0,a_1\}),\cone(\{a_1,a_2\}),\cone(\{a_2,a_0\})\}
\end{align*}
\item[(ii)] The fans $F_a$ representing the Hirzebruch surfaces for $a\in\Z$.
\begin{align*}
F_a^{(1)}&= \{a_0,a_1,a_2,a_3\} = \{(1,0),(0,1),(-1,-a),(0,-1)\} \\
F_a^{(d)}&=\{\cone(\{a_0,a_1\}),\cone(\{a_1,a_2\}),\cone(\{a_2,a_3\}),\cone(\{a_3,a_0\})\}
\end{align*}
\end{itemize}
\end{theorem}
\begin{center}
\begin{tikzpicture}[scale=0.5]
\draw[black!70,->] (0,0) -- (3,0);
\draw[black!70,->] (0,0) -- (0,3);
\draw[black!70] (0,0) -- (-3,0);
\draw[black!70] (0,0) -- (0,-3);
\draw[black,thick,->] (0,0) -- (2,0) node[anchor=south] {\tiny $(1,0)$};
\draw[black,thick,->] (0,0) -- (0,2) node[anchor=west] {\tiny $(0,1)$};
\draw[black,thick,->] (0,0) -- (-2,-2) node[anchor=north west] {\tiny $(-1,-1)$};
\draw (0,0) node[anchor=south west] {$F_p$};
\end{tikzpicture}
\quad\quad
\begin{tikzpicture}[scale=0.5]
\draw[black!70,->] (0,0) -- (3,0);
\draw[black!70,->] (0,0) -- (0,3);
\draw[black!70] (0,0) -- (-3,0);
\draw[black!70] (0,0) -- (0,-3);
\draw[black,thick,->] (0,0) -- (2,0) node[anchor=south] {\tiny $(1,0)$};
\draw[black,thick,->] (0,0) -- (0,2) node[anchor=west] {\tiny $(0,1)$};
\draw[black,thick,->] (0,0) -- (0,-2) node[anchor=west] {\tiny $(0,-1)$};
\draw[black,thick,->] (0,0) -- (-2,-1) node[anchor=north] {\tiny $(-1,-a)$};
\draw (0,0) node[anchor=south west] {$F_a$};
\end{tikzpicture}
\end{center}

\paragraph{Tree of Blowing-ups}
Instead of generating combinatorial types or finding parameters we have a different task now, namely equivariant blowing-ups. In dimension two this means choosing one full-dimensional cone $C=\cone(r_1,r_2)$, inserting the ray $r=r_1+r_2$ and replacing the cone $C$ with $C_1=\{r_1,r\}$ and $C_2=\{r,r_2\}$. \\
So beginning with a fan $F$ we can generate all normal fans with the following recursion:
\begin{enumerate}
\item For every full-dimensional cone $C\in F^{(d)}$:
\begin{enumerate}
\item Create the blowing-up $b(F,C)$,
\item If $b(F,C)$ has less than $\nlimit$ full-dimensional cones: Start recursion.
\end{enumerate}
\end{enumerate}
Which is a depth-first search in a rather large rooted tree, having $1 + (3!+4!+\cdots+11!)/2! = 21977356$ nodes for up to $12$ full-dimensional cones for the fan $F_p$, but lots or even most of those fans are isomorphic. So we need a way to reduce the number of nodes to consider:
\par
We first assign an order to all full-dimensional cones. Now, since blowing up $b(F,C_i)$ in $C_j$ and $b(F,C_j)$ in $C_i$ will give the same fan, we can ignore all cones $C_j$ with index $j<i$ in the subtree of the blowing-up $b(F,C_i)$. (When blowing up in cone $C_i$ in a fan with $k$ full-dimensional cones $C_1,\ldots,C_k$ we will index the two new cones with $i$ and $k+1$.)\\
This reduces the amount of blowing-ups to $58785$ for $F_p$ and $35072$ for $F_a$.

\paragraph{Parameters}
The fan $F_a$ and all further blowing-ups of it represent classes of infintely many fans with $a\in\Z$. But, from the proof of our main theorem, we know that for some upper bound $\limit$ on the number of lattice points, there can be only finitely many of those. The proof even gives us a direct way to bound the parameter $a$.
\begin{lem}
For the fan $F_a$ and a maximal number of lattice points $\limit$ for the polytope, we can restrict $a$ to $0 \le a \le \limit$ and $a \ne 1$.
\end{lem}
\begin{pf}
For a pair of two-dimensional cones $C_1=\{r_1,r\}$ and $C_2=\{r_2,r\}$ sharing a common ray $r$ we can write $r_1+r_2$ in terms of $r$, i.e.~$r_1+r_2 = a_r r$. Additionally we know $a_r\ge-l(e)$, where $l(e)$ is the length of the edge $e$ dual to the ray $r$. Obviously $l(e)\le\limit$, so we have $a_r\ge-\limit$. 
\begin{itemize}
\item For the cones $C_1=\{(1,0),(0,1)\}$, $C_2=\{(-1,-a),(0,1)\}$ this yields
\[
r_1+r_2=(1,0)+(-1,-a)=(0,-a)=-a\cdot(0,1)
\]
and thus $-a = a_{(0,1),1}\ge-\limit$.
\item And for $C_1=\{(1,0),(0,-1)\}$, $C_2=\{(-1,-a),(0,-1)\}$ we have 
\[
r_1+r_2=(1,0)+(-1,-a)=(0,-a)=a\cdot(0,-1)
\] and $a = a_{(0,-1),1}\ge-\limit$.
\end{itemize} 
Because of the symmetry of the fan we can restrict to positive $a$ and skip $a=1$ since $F_a$ then equals $F_p$ after an equivariant blowing-up in the cone $C = \cone(\{a_0,a_2\})$.
\qed
\end{pf}

\paragraph{Parameters in the Blowing-up}
We could find bounds on the parameter $a$ of the Hirzebruch surface using an edge of the corresponding polytope, but when blowing up the fan, which means cutting off a vertex of the polytope, this edge may have different parameters and the bounds may be different too. We could of course do the same calculation for all blowing-ups, but a faster approach is to find a way to use the bounds we found for some fan $F$ on all its blowing-ups.\\
Assume we have a pair of cones $C_1=\{r_1,r\},\ C_2=\{r_2,r\}\in F$ around a ray $r$ with $r_1+r_2 = a_r r$, this means our bound is $a_r\ge-\limit$. Now, in the blowing-up $b(F,C_2)$, we have a new cone $C_2'=\{r_2',r\}$ with $r_2' = r+r_2$. This yields 
\[
r_1+r_2' = a_r r - r_2 + r + r_2 = (a_r+1)r
\]
So we now have $a_r+1\ge-\limit$ and thus $a_r\ge-\limit-1$. This means we have to loosen the bounds on a parameter by one if it appears as edge-parameter at one of the rays of the current cone. \\
Thus, after creating a new fan for every valid value of $a$ for each class of fans, we now have all possible smooth normal fans for the smooth polytopes.

\subsubsection{Creating the Polytopes}
Now we can continue as in the general case, for every fan $F$ we define the polytope $\B(F,N)$ approximating $\RHS(F,N)$ and enumerate all its lattice points $b$. For some fan with rays $A$ and a right-hand side $b$ check if $|P(A,b)\cap\Z^2| \le \limit$. And finally remove polytopes that are lattice-isomorphic to another one.

\subsubsection{Results}
We can now state the result of the computation:
\begin{prop}
There are 41 smooth lattice polygons with at most \nlimit lattice points.
\begin{center}
\begin{tabular}{r|cccccccccc}
\toprule
Vertices 		& 3 &  4 & 5 & 6 & 7 & 8 & $\ge$9 \\
\hline
Polygons 			& 3 & 30 & 3 & 4 & 0 & 1 &    0 \\
\bottomrule
\end{tabular}
\end{center}
The complete list of polytopes can be found in the appendix.
\end{prop}

\subsection{Smooth Polytopes in Dimension Three}
\subsubsection{Classifying the Normal Fans}
\paragraph{Minimal Fans}
As in the two-dimensional case, we start with a list of smooth minimal fans from Tadao Oda~\cite{oda}.

\begin{theorem}{\bf(Tadao Oda~\cite{oda})}
Every smooth fan in dimension three, with at most eight rays, can, up to isomorphism, be created by finitely many successive equivariant blowing-ups of one of the $19$ fans with labels $3^4$, $(3^24^3)'$, $(3^24^3)''$, $4^6$, $3^24^36^2$, $3^14^35^3$, $4^55^2$, $4^66^2$, $3^24^47^2$, $3^34^15^16^3$, $3^24^25^26^2$, $3^14^45^16^2$, $3^24^15^46^1$, $(3^14^35^36^1)'$, $(3^14^35^36^1)''$, $(3^25^6)'$, $(3^25^6)''$, $(4^45^4)'$ and $(4^45^4)''$ and weights as given in \emph{Convex Bodies and Algebraic Geometry}.
\end{theorem}

\begin{rem}\textcolor{white}{i}
\begin{itemize}
\item All drawings of such fans in the following chapter are stereographic projections of the corresponding triangulation of $S^2$.
\item The labels give the occurence of vertices with different degrees in the triangulation, for example $3^24^36^2$ means $2$ vertices of degree $3$, $3$ vertices of degree $4$ and $2$ vertices of degree $6$. See the pictures later in this chapter.
\end{itemize}
\end{rem}
\begin{rem} With at most eight rays in the fan, from Eulers formula and the fact the the fan is simplicial we know that there are at most twelve full-dimensional cones in the fan and thus at most twelve vertices in the polytope. This is the main reason we chose $12$ as upper bound for the number of lattice points in the polytope.
\end{rem}

Even in the two-dimensional case, the blowing-up tree was huge so it comes in handy that we can use the two-dimensional classification now to reduce the number of fans we have to consider. From the $19$ fans of the theorem there will be only five left after applying the now following criterion and it will help us a lot later during the computation of all blowing-ups.

\paragraph{Smooth Polygons}
We use the smooth polygons we generated earlier to remove fans whose polytopes would have too many lattice points on the boundary. Let $\mathcal{P}_N(k)$ be the set of all smooth lattice $k$-gons with at most $N$ lattice points. Define the following numbers:
\begin{align*}
l_N(k) &= \min\{|P\cap\Z^2|\ \mbox{for} \ P\in\mathcal{P}_N(k)\} \\
i_N(k) &= \min\{|\inter(P)\cap\Z^2|\ \mbox{for} \ P\in\mathcal{P}_N(k)\} \\
b_N(k) &= \min\{|P\cap\Z^2| - |\inter(P)\cap\Z^2|\ \mbox{for} \ P\in\mathcal{P}_N(k)\} 
\end{align*}
From the results of the previous section we can determine the values shown in the table below.
\begin{center}
\begin{tabular}{r|cccccccccc}
\toprule
$k$ 			& 3 & 4 & 5 & 6 &   7 &  8 &   $>$9 \\
\hline
$l_{12}(k)$ 	& 3 & 4 & 8 & 7 & $>$12 & 12 & $>$12  \\
$i_{12}(k)$	& 0 & 0 & 1 & 1 &     &  4 &     \\
$b_{12}(k)$	& 3 & 4 & 7 & 6 &     &  8 &     \\
\bottomrule
\end{tabular}
\end{center}
Denote by $k_r$ the number of full-dimensional cones in $F$ a ray $r$ is contained in and let $k_F$ be the maximal $k_r$ over all rays $r$ in $F$. \\
All vertices of the polytope defined by the full-dimensional cones around a ray $r$ lie on the boundary of the facet defined by the ray $r$, so for every ray $r$ there is a $k_r$-gon on the boundary of the polytope. Since the polytope is smooth, all these polygons have to be smooth, too. This yields the following lower bound for the number of lattice points, which depends solely on the combinatorial type of the fan. 

\begin{lem}{\bf (smooth lattice polygon criterion)}
For a smooth lattice polytope $P$ with the normal fan $F$.
\[ 
|P\cap\Z^3| \ge |F^{(d)}| + (b_N(k_F) - k_F) + \sum_{r\in F^{(1)}} i_N(k_r)
\]
\end{lem}
\begin{rem}
This bound can be improved by considering the combinatorial type of the fan more precisely. However, for our computation it suffices.
\end{rem}

\paragraph{Remaining Minimal Fans}
Now we can check the fans from the theorem for the above criterion, the labels of the fans are exactly the numbers $k_r^*$ we are looking for.
\begin{asparaitem}
\item First, consider all fans with eight rays and twelve full-dimensional cones. That are the fans labeled $4^66^2$, $3^24^47^2$, $3^34^15^16^3$, $3^24^25^26^2$, $3^14^45^16^2$, $3^24^15^46^1$, $(3^14^35^36^1)'$, $(3^14^35^36^1)''$, $(3^25^6)'$, $(3^25^6)''$, $(4^45^4)'$ and $(4^45^4)''$, which cannot be blown up as this would add two additional full-dimensional cones. From the labels one can see that they all have at least one ray contained in five, six or seven full-dimensional cones ($5^*,6^*$ or $7^*$), this means the corresponding polytope would have a five-, six- or seven-gon on the boundary. But since any smooth five- or six-gon has at least one interior lattice point and there is no smooth seven-gon with less than $13$ lattice points at all, the polytope would have at least $13$ lattice points on the boundary. So none of these fans admits a smooth lattice polytope with at most $\nlimit$ lattice points.
\item Now consider the fans labeled $4^55^2$ and $3^14^35^3$ with ten full-dimensional cones.
\begin{itemize}
\item[$(4^55^2)$] 
\parbox[t]{0.97\linewidth}{
\parbox[t]{0.81\linewidth}{
Here we have two non-incident rays each contained in five full-dimensional cones, which means we have two edge-disjoint five-gons in the polytope. But we know that there is no smooth five-gon with less than eight lattice points and thus any smooth polytope with this fan as normal fan would have at least $16$ vertices.}
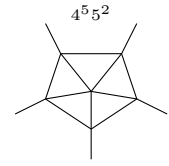
\parbox[t]{0.16\linewidth}{
\begin{center}
{\tiny $4^55^2$}\\
\begin{tikzpicture}[scale=0.2]
\draw (3.00,7.00)--(7.00,7.00)--(8.00,4.00)--(5.00,2.00)--(2.00,4.00)--cycle;
\draw (5.00,2.00)--(5.00,4.50);
\draw (2.00,4.00)--(5.00,4.50);
\draw (3.00,7.00)--(5.00,4.50);
\draw (7.00,7.00)--(5.00,4.50);
\draw (8.00,4.00)--(5.00,4.50);
\draw (5.00,2.00)--(5.00,0.00);
\draw (8.00,4.00)--(10.00,3.00);
\draw (2.00,4.00)--(0.00,3.00);
\draw (3.00,7.00)--(2.00,9.00);
\draw (7.00,7.00)--(8.00,9.00);
\end{tikzpicture}\\
\end{center}
}\\
A smooth polytope corresponding to any blowing-up of this fan would have $12$ vertices and either two six-gons or a five-gon and a six-gon and thus have at least $14$ lattice points.
}
\item[$(3^14^35^3)$] \parbox[t]{0.76\linewidth}{
This fan contains three rays with five neighbors, so any polytope with this normal fan would have at least three lattice points in the relative interior of the facets and thus at least $13$ lattice points in total. Any blowing-up would leave at least one of the five-gons unchanged and thus there would still be at least $13$ lattice points, since we have now two additional vertices from the blowing-up and one in the interior lattice points in the five-gon.}
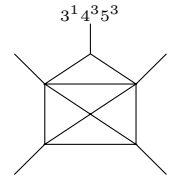
\parbox[t]{0.2\linewidth}{
\begin{center}
{\tiny $3^14^35^3$}\\
\begin{tikzpicture}[scale=0.2]
\draw (0.00,0.00)--(2.00,2.00)--(2.00,6.00)--(5.00,4.00)--(2.00,2.00)--
(8.00,2.00)--(5.00,4.00)--(8.00,6.00)--(2.00,6.00)--(5.00,8.00)--(8.00,6.00)--
(8.00,2.00)--(10.00,0.00);
\draw (2.00,6.00)--(0.00,8.00);
\draw (8.00,6.00)--(10.00,8.00);
\draw (5.00,8.00)--(5.00,10.00);
\end{tikzpicture}
\end{center}}
\end{itemize}
\item By the above criterion, the fans numbered $3^4, (3^24^3)', (3^24^3)'', 4^6$ and $3^24^36^2$ could appear as normal fans for smooth lattice polytopes with at most $\nlimit$ lattice points and will be the input to the algorithm.
\end{asparaitem}
This leaves us with the following five fans or classes thereof with rays and edge-para\-meters as drawn~\cite{oda}.\\
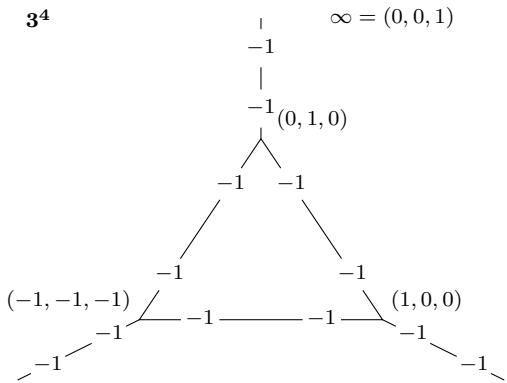
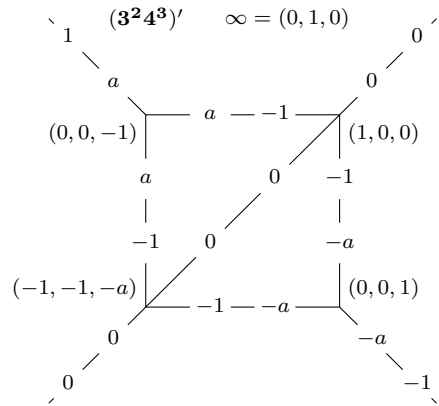
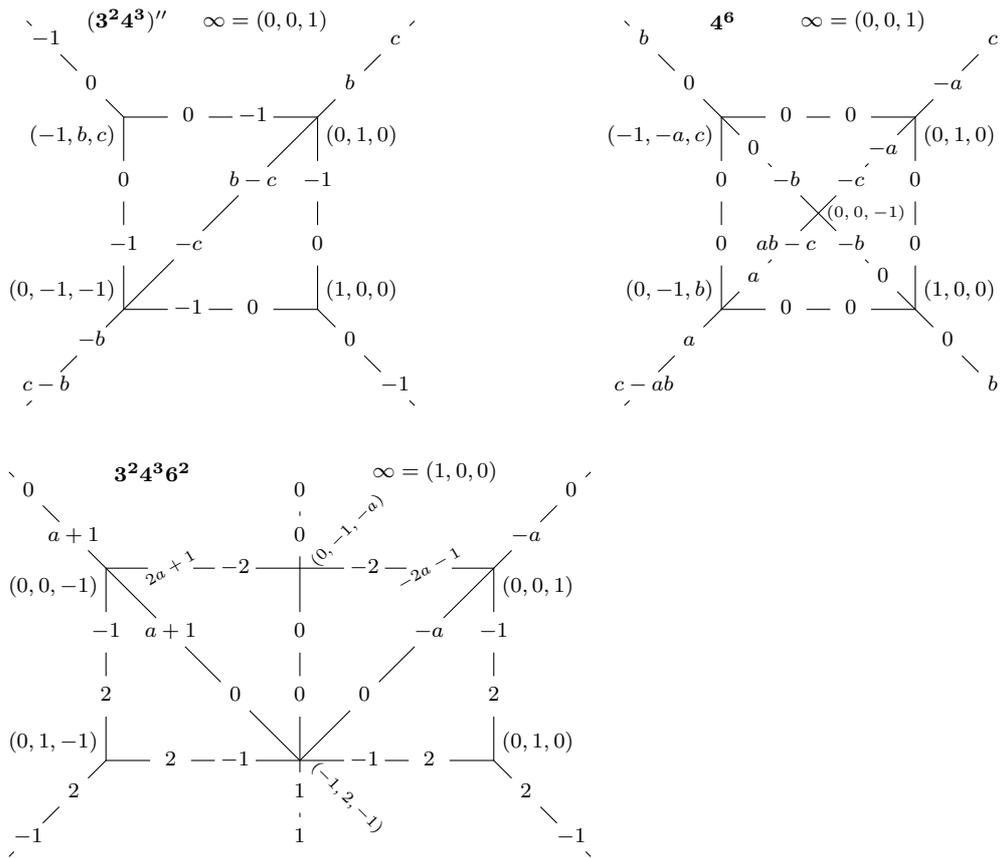
\begin{center}
{\scriptsize
\begin{tabularx}{\linewidth}{XX}
\begin{tikzpicture}[scale=0.8]
\draw (2.00,1.00)--(6.00,1.00)--(4.00,4.00)--cycle;
\draw (4.00,6.00)--(4.00,4.00) node[anchor=south west] {$\ (0,1,0)$};
\draw (8.00,0.00)--(6.00,1.00) node[anchor=south west] {$(1,0,0)$};
\draw (0.00,0.00)--(2.00,1.00) node[anchor=south east] {$(-1,-1,-1)$};
\draw (0,6) node[anchor=west] {$\mathbf{3^4}$};
\draw (5,6) node[anchor=west] {$\infty=(0,0,1)$};
\edgeparam{3,1}{-1}
\edgeparam{5,1}{-1}
\edgeparam{5.5,1.75}{-1}
\edgeparam{4.5,3.25}{-1}
\edgeparam{2.5,1.75}{-1}
\edgeparam{3.5,3.25}{-1}
\edgeparam{4,4.5}{-1}
\edgeparam{4,5.5}{-1}
\edgeparam{7.5,0.25}{-1}
\edgeparam{6.5,0.75}{-1}
\edgeparam{0.5,0.25}{-1}
\edgeparam{1.5,0.75}{-1}

\end{tikzpicture} &
\begin{tikzpicture}[scale=0.85]

\draw (1.5,1.5) -- (4.5,1.5) -- (4.5,4.5) -- (1.5,4.5) -- cycle;
\draw (0,0) -- (1.5,1.5) node[anchor=south east] {$(-1,-1,-a)$};
\draw (6,0) -- (4.5,1.5) node[anchor=south west] {$(0,0,1)$};
\draw (6,6) -- (4.5,4.5) node[anchor=north west] {$(1,0,0)$};
\draw (0,6) -- (1.5,4.5) node[anchor=north east] {$(0,0,-1)$};
\draw (1.5,1.5) -- (4.5,4.5);
\draw (0.8,6) node[anchor=west] {$\mathbf{(3^24^3)'}$};
\draw (2.6,6) node[anchor=west] {$\infty=(0,1,0)$};

\edgeparam{0.3,0.3}{0}
\edgeparam{1,1}{0}
\edgeparam{2.5,2.5}{0}
\edgeparam{3.5,3.5}{0}
\edgeparam{5,5}{0}
\edgeparam{5.7,5.7}{0}
\edgeparam{0.3,5.7}{1}
\edgeparam{1,5}{a}
\edgeparam{1.5,2.5}{-1}
\edgeparam{1.5,3.5}{a}
\edgeparam{2.5,4.5}{a}
\edgeparam{3.5,4.5}{-1}
\edgeparam{5.7,0.3}{-1}
\edgeparam{5,1}{-a}
\edgeparam{2.5,1.5}{-1}
\edgeparam{3.5,1.5}{-a}
\edgeparam{4.5,2.5}{-a}
\edgeparam{4.5,3.5}{-1}

\end{tikzpicture}
\end{tabularx}
\\ \textcolor{white}{x} \\ \textcolor{white}{x} \\
\begin{tabularx}{\linewidth}{XX}
\begin{tikzpicture}[scale=0.85]
\draw (1.5,1.5) -- (4.5,1.5) -- (4.5,4.5) -- (1.5,4.5) -- cycle;
\draw (0,0) -- (1.5,1.5) node[anchor=south east] {$(0,-1,-1)$};
\draw (6,0) -- (4.5,1.5) node[anchor=south west] {$(1,0,0)$};
\draw (6,6) -- (4.5,4.5) node[anchor=north west] {$(0,1,0)$};
\draw (0,6) -- (1.5,4.5) node[anchor=north east] {$(-1,b,c)$};
\draw (1.5,1.5) -- (4.5,4.5);
\draw (0.8,6) node[anchor=west] {$\mathbf{(3^24^3)''}$};
\draw (2.6,6) node[anchor=west] {$\infty=(0,0,1)$};

\edgeparam{0.3,0.3}{c-b}
\edgeparam{1,1}{-b}
\edgeparam{2.5,2.5}{-c}
\edgeparam{3.5,3.5}{b-c}
\edgeparam{5,5}{b}
\edgeparam{5.7,5.7}{c}
\edgeparam{0.3,5.7}{-1}
\edgeparam{1,5}{0}
\edgeparam{1.5,2.5}{-1}
\edgeparam{1.5,3.5}{0}
\edgeparam{2.5,4.5}{0}
\edgeparam{3.5,4.5}{-1}
\edgeparam{5.7,0.3}{-1}
\edgeparam{5,1}{0}
\edgeparam{2.5,1.5}{-1}
\edgeparam{3.5,1.5}{0}
\edgeparam{4.5,2.5}{0}
\edgeparam{4.5,3.5}{-1}

\end{tikzpicture} &
\begin{tikzpicture}[scale=0.85]
\draw (1.5,1.5) -- (4.5,1.5) -- (4.5,4.5) -- (1.5,4.5) -- cycle;
\draw (0,0) -- (1.5,1.5) node[anchor=south east] {$(0,-1,b)$};
\draw (6,0) -- (4.5,1.5) node[anchor=south west] {$(1,0,0)$};
\draw (6,6) -- (4.5,4.5) node[anchor=north west] {$(0,1,0)$};
\draw (0,6) -- (1.5,4.5) node[anchor=north east] {$(-1,-a,c)$};
\draw (1.5,1.5) -- (4.5,4.5);
\draw (1.5,4.5) -- (4.5,1.5);
\draw (3,3) node[anchor=west] {\tiny$(0,0,-1)$};
\draw (1.2,6) node[anchor=west] {$\mathbf{4^6}$};
\draw (2.6,6) node[anchor=west] {$\infty=(0,0,1)$};

\edgeparam{0.3,0.3}{c-ab}
\edgeparam{1,1}{a}
\edgeparam{2,2}{a}
\edgeparam{2.5,2.5}{ab-c}
\edgeparam{3.5,3.5}{-c}
\edgeparam{4,4}{-a}
\edgeparam{5,5}{-a}
\edgeparam{5.7,5.7}{c}
\edgeparam{2,4}{0}
\edgeparam{2.5,3.5}{-b}
\edgeparam{3.5,2.5}{-b}
\edgeparam{4,2}{0}

\edgeparam{0.3,5.7}{b}
\edgeparam{1,5}{0}
\edgeparam{1.5,2.5}{0}
\edgeparam{1.5,3.5}{0}
\edgeparam{2.5,4.5}{0}
\edgeparam{3.5,4.5}{0}
\edgeparam{5.7,0.3}{b}
\edgeparam{5,1}{0}
\edgeparam{2.5,1.5}{0}
\edgeparam{3.5,1.5}{0}
\edgeparam{4.5,2.5}{0}
\edgeparam{4.5,3.5}{0}
\end{tikzpicture} 
\end{tabularx}
\\ \textcolor{white}{x} \\ \textcolor{white}{x} \\
\begin{tabularx}{\linewidth}{X}
\begin{tikzpicture}[scale=0.85]
\draw (1.5,1.5) -- (7.5,1.5) -- (7.5,4.5) -- (1.5,4.5) -- cycle;
\draw (0,0) -- (1.5,1.5) node[anchor=south east] {$(0,1,-1)$};
\draw (9,0) -- (7.5,1.5) node[anchor=south west] {$(0,1,0)$};
\draw (9,6) -- (7.5,4.5) node[anchor=north west] {$(0,0,1)$};
\draw (0,6) -- (1.5,4.5) node[anchor=north east] {$(0,0,-1)$};
\draw (4.5,1.5) -- (7.5,4.5);
\draw (1.5,4.5) -- (4.5,1.5);
\draw (4.5,0) -- (4.5,6);
\draw (1.5,6) node[anchor=west] {$\mathbf{3^24^36^2}$};
\draw (5.5,6) node[anchor=west] {$\infty=(1,0,0)$};

\edgeparam{0.3,0.3}{-1}
\edgeparam{1,1}{2}
\edgeparam{5.5,2.5}{0}
\edgeparam{6.5,3.5}{-a}
\edgeparam{8,5}{-a}
\edgeparam{8.7,5.7}{0}
\edgeparam{0.3,5.7}{0}
\edgeparam{1,5}{a+1}
\edgeparam{2.5,3.5}{a+1}
\edgeparam{3.5,2.5}{0}
\edgeparam{8.7,0.3}{-1}
\edgeparam{8,1}{2}

\edgeparam{1.5,2.5}{2}
\edgeparam{1.5,3.5}{-1}
\edgeparam{4.5,0.3}{1}
\edgeparam{4.5,1}{1}
\edgeparam{4.5,2.5}{0}
\edgeparam{4.5,3.5}{0}
\edgeparam{4.5,5}{0}
\edgeparam{4.5,5.7}{0}
\edgeparam{7.5,2.5}{2}
\edgeparam{7.5,3.5}{-1}
\edgeparam{2.5,1.5}{2}
\edgeparam{3.5,1.5}{-1}
\edgeparam{5.5,1.5}{-1}
\edgeparam{6.5,1.5}{2}
\edgeparamx{2.5,4.5}{2a+1}
\edgeparam{3.5,4.5}{-2}
\edgeparam{5.5,4.5}{-2}
\edgeparamx{6.5,4.5}{-2a-1}

\draw (4.8,1.7) node[rotate=-45,anchor=north west] {\tiny$(-1,2,-1)$};
\draw (4.8,4.3) node[rotate=45,anchor=south west] {\tiny$(0,-1,-a)$};
\end{tikzpicture}
\end{tabularx}
}
\end{center}

\begin{figure}[h]
\begin{center}
\begin{tikzpicture}[scale=1.0]
\draw[black!60] (1,0) -- (3.7,0);
\draw[black!40] (1,0) -- (4.5,0.875);
\draw[black!60] (1,0) -- (4,2.5);
\draw[black!60] (2.4,0) -- (3,1.67) -- (4,0.75) -- cycle;
\draw[black,dashed,thick] (3,1.67) -- (3.2,0.9);
\draw[black,dashed,thick] (2.4,0) -- (3.2,0.9);
\draw[black,dashed,thick] (4,0.75) -- (3.2,0.9);
\draw[black,dashed,thick] (1,0) -- (4.5,1.43);
\fill[black] (3.2,0.9) circle (2pt);

\draw[<->,black] (5,1) -- (5.5,1);
\end{tikzpicture}
\begin{tikzpicture}[scale=1.0]
\draw[black!60] (1.5,1.5) -- (1.5,0);
\draw[black!60] (1.5,1.5) -- (0,2.2);
\draw[black,dashed,thick] (1.48,0.6) -- (0.8,1.85) -- (0.6,1.1) -- cycle;
\draw[black!60] (1.5,1.5) -- (0,0.8);
\draw[white] (2,2) circle (2pt);
\end{tikzpicture}
\quad\quad
\begin{tikzpicture}[scale=1.0]
\draw[black!40] (1,0) -- (4.5,0.875);
\draw[black!60] (1,0) -- (4,2.5);
\draw[black!60] (1,0) -- (4,-0.5);
\draw[black!60] (2.4,0) -- (3,1.67) -- (4,0.75) -- 
					 (3.3,-0.38) -- (2.4,0) -- (4,0.75);
\draw[black,dashed,thick] (3,1.67) -- (3.3,0.42);
\draw[black,dashed,thick] (2.4,0) -- (3.3,0.42);
\draw[black,dashed,thick] (4,0.75) -- (3.3,0.42);
\draw[black,dashed,thick] (3.3,-0.38) -- (3.3,0.42);
\draw[black,dashed,thick] (1,0) -- (4.5,0.64);
\fill[black] (3.3,0.42) circle (2pt);
\draw[black!60] (1,0) -- (3.7,0);

\draw[<->,black] (5,0.5) -- (5.5,0.5);
\end{tikzpicture}
\begin{tikzpicture}[scale=1.0]
\draw[black!40] (0,0) -- (1.5,0.3);
\draw[black!40] (0.6,2) -- (1.5,1.1);
\draw[black!60] (0.3,-0.3) -- (1.5,0.3) -- (1.5,1.1) -- (0.3,1.7);
\draw[black,dashed,thick]  (0.8,-.05) -- (0.8,1.45) -- (1.0,1.6) -- 
									(1.0,0.2) -- cycle;
\end{tikzpicture}
\caption{Blowing-up of a three-dimensional fan and the cutting-off in the corresponding polytope}
\label{cutoff}
\end{center}
\end{figure}
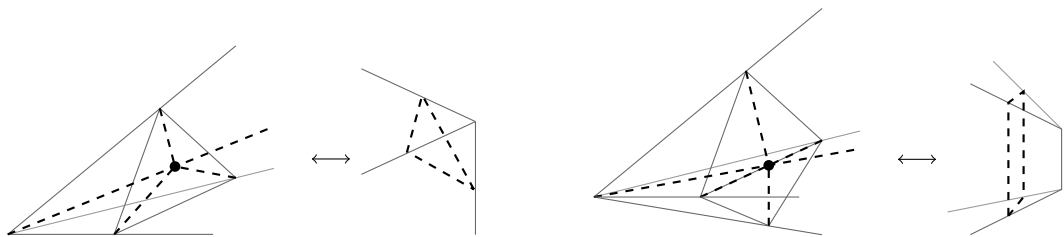

\paragraph{Tree of Blowing-ups}
Again, the first step of the algorithm is to generate all possible normal fans by successive blowing-ups. There are now two possible cases, blowing up in a $3$-dimensional cone or a $2$-dimensional cone, which corresponds in the polytope to cutting off a vertex or an edge respectively, as shown in figure~\ref{cutoff}.

Generating all blowing-ups successively gives again a huge rooted tree, but we now have more ways to reduce the number of nodes we have to check.
\begin{itemize}
\item[(1)] We can again limit the depth of the tree since we want at most $\nlimit$ full-dimensional cones and each blowing-up increases the number of such cones by two.
\item[(2)] Blowing-ups in full-dimensional cones can again be done in any order. Let $F$ contain $k$ full-dimensional cones. When blowing up in the cone $C_i$ the three new cones get the numbers $i,k+1,k+2$, so we can ignore all full-dimensional cones with index strictly less than $i$ at all nodes in the subtree of this blowing-up.
\item[(3)] For codimension $1$ cones there are two cases.
Let $C$ be a codimension $1$ cone contained in two full-dimensional cones $C_1$ and $C_2$.
\begin{itemize}
\item[(3a)] Blowing up $b(F,C)$ in another codimension $1$ cone $C'$, which is contained in $C_1$ or $C_2$, will yield a different result than blowing up $b(F,C')$ in the cone $C$. All cones have to be considered in this case.
\item[(3b)] Blowing up $b(F,C)$ in a codimension $1$ cone $C'$, which is not contained in $C_1$ or $C_2$, will yield the same result as blowing up $b(F,C')$ in $C$. So we can use an ordering on the codimension $1$ cones to skip cones with lower index if they do not fall into case $(3a)$.
\end{itemize}
\item[(4)] Blowing up in a full-dimensional cone and later in one of its codimension $1$ faces yields the same rays and the same combinatorial type as blowing up in the codimension $1$ cone first and then in on new codimension $1$ cone as shown below.
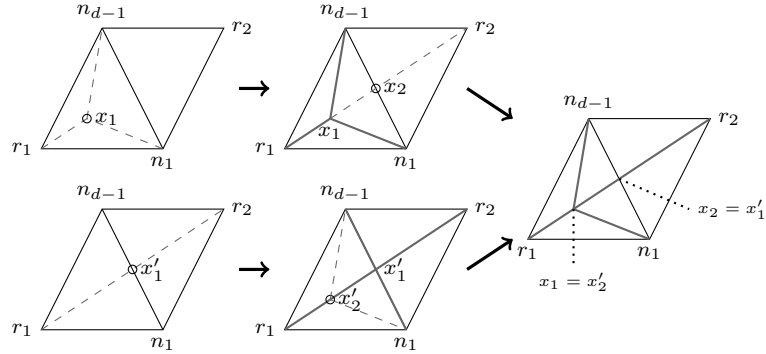
\begin{center}
{ \scriptsize
\begin{tikzpicture}[scale=0.4]
\draw[black] (-8,0) -- (-4,0) -- (-2,4) -- (-6,4) -- cycle;
\draw[black] (-4,0) -- (-6,4);
\draw[black!60,dashed] (-8,0) -- (-2,4);
\draw (-5,2) circle (4pt) node[anchor=west] {$x_1'$};
\draw (-8,0) node[anchor=east] {$r_1$};
\draw (-2,4) node[anchor=west] {$r_2$};
\draw (-6,4) node[anchor=south] {$n_{d-1}$};
\draw (-4,0) node[anchor=north] {$n_1$};

\draw[black] (-8,6) -- (-4,6) -- (-2,10) -- (-6,10) -- cycle;
\draw[black] (-4,6) -- (-6,10);
\draw[black!60,dashed] (-8,6) -- (-6.5,7) -- (-4,6);
\draw[black!60,dashed] (-6.5,7) -- (-6,10);
\draw (-6.5,7) circle (4pt) node[anchor=west] {$x_1$};
\draw (-8,6) node[anchor=east] {$r_1$};
\draw (-2,10) node[anchor=west] {$r_2$};
\draw (-6,10) node[anchor=south] {$n_{d-1}$};
\draw (-4,6) node[anchor=north] {$n_1$};

\draw[very thick,->] (-1.5,2) -- (-0.5,2);
\draw[very thick,->] (-1.5,8) -- (-0.5,8);

\draw[black] (0,6) -- (4,6) -- (6,10) -- (2,10) -- cycle;
\draw[black] (4,6) -- (2,10);
\draw[black!60,thick] (0,6) -- (1.5,7) -- (2,10);
\draw[black!60,thick] (1.5,7) -- (4,6);
\draw (3,8) circle (4pt);
\draw[black!60,dashed] (1.5,7) -- (6,10);
\draw (0,6) node[anchor=east] {$r_1$};
\draw (6,10) node[anchor=west] {$r_2$};
\draw (2,10) node[anchor=south] {$n_{d-1}$};
\draw (4,6) node[anchor=north] {$n_1$};
\draw (1.5,7) node[anchor=north] {$x_1$};
\draw (3,8) node[anchor=west] {$x_2$};
\draw[black] (0,0) -- (4,0) -- (6,4) -- (2,4) -- cycle;
\draw[black] (4,0) -- (2,4);
\draw[black!60,thick] (0,0) -- (6,4);
\draw[black!60,thick] (2,4) -- (4,0);
\draw (1.5,1) circle (4pt);
\draw[black!60,dashed] (1.5,1) -- (2,4);
\draw[black!60,dashed] (1.5,1) -- (4,0);
\draw (0,0) node[anchor=east] {$r_1$};
\draw (6,4) node[anchor=west] {$r_2$};
\draw (2,4) node[anchor=south] {$n_{d-1}$};
\draw (4,0) node[anchor=north] {$n_1$};
\draw (1.5,1) node[anchor=west] {$x_2'$};
\draw (3,2) node[anchor=west] {$x_1'$};

\draw[very thick,->] (6,8) -- (7.5,7);
\draw[very thick,->] (6,2) -- (7.5,3);

\draw[black] (8,3) -- (12,3) -- (14,7) -- (10,7) -- cycle;
\draw[black] (12,3) -- (10,7);
\draw[black!60,thick] (8,3) -- (14,7);
\draw[black!60,thick] (12,3) -- (9.5,4) -- (10,7);
\draw (8,3) node[anchor=north] {$r_1$};
\draw (14,7) node[anchor=west] {$r_2$};
\draw (10,7) node[anchor=south] {$n_{d-1}$};
\draw (12,3) node[anchor=north] {$n_1$};

\draw[thick,dotted] (9.5,4) -- (9.5,2.2) node[anchor=north] {\tiny $x_1=x_2'$};
\draw[thick,dotted] (11,5) -- (13.3,4) node[anchor=west] {\tiny $x_2=x_1'$};

\end{tikzpicture} }
\end{center}
And we have the following new rays:
\begin{align*}
x_2&=\sum_{i=1}^{d-1}n_i 	  & x_1' &= \sum_{i=1}^{d-1}n_i = x_2 \\
x_1&=r_1+\sum_{i=1}^{d-1}n_i & x_2' &= x_1'+r_1=r_1+\sum_{i=1}^{d-1}n_i = x_1 
\end{align*} So only one of these cases has to be considered.
\item[(5)] As already mentioned, we can use the smooth lattice polygon criterion described earlier. Fans that violate this criterion will be considered for further blowing-ups, but will not be kept for the later steps of the algorithm. We cannot cut off the tree here because the criterion is not monotone:\\
Assume we have a polytope with one vertex that is contained in three five-gons, then cutting-off this vertex will replace all with six-gons. So in total this new polytope will have two more vertices but may have less lattice points than the original polytope, since the smallest smooth five-gon has $8$ lattice points in total and the smallest smooth six-gon is the hexagon with only $7$ total lattice points.
\end{itemize}
This has been implemented with one variable for each full-dimensional and each codimension $1$ cone that can have the values `always ignore', `ignore', `consider', `always consider'. We will create the blowing-up and continue the recursion only for cones with value `consider' or `always consider'. The specific cases are:
\begin{itemize}
\item[(2)] Visited full-dimensional cones are marked `ignore' for the subtree.
\item[(3a)] Visited codimension $1$ cones are marked `ignore' for the subtree if they haven't been marked `always consider' or `always ignore'.
\item[(3b)] Codimension $1$ cones are marked `always consider' if they were on `ignore' before and are contained in one of the full-dimensional cones containing the codimension $1$ cone of the current blowing-up.
\item[(4)] Codimension $1$ cones are marked `always ignore' if they were contained in the full-dimensional cone of the current blowing-up.
\item[(x)] Newly created cones are always set to `consider'.
\end{itemize}

\paragraph{Parameters for the Minimal Fans}
Most of the remaining minimal fans are actually classes of fans with parameters $a,b,c\in\Z$. Hence we have to find bounds for their parameters and find ways to use those bounds for their blowing-ups.\\
We can find bounds for their parameters analogously to the two-dimensional case: Whenever a parameter $x$ appears as $a_{C,i}=x$ and $a_{C',j}=-x$ we know $-\limit \le x \le \limit $ and get:
\begin{align*}
[3^4] &\quad \mbox{has no parameters,}\\
[(3^24^3)'] & \quad 0\le a \le \limit \mbox{ using the symmetry  of the fan,} \\
[(3^24^3)''] & \quad -\limit\le b,c \le \limit, \\
[4^6] & \quad -\limit\le a,b,c \le \limit, \\
[3^24^36^2] & \quad \mbox{here we use $a_{C,i}=2a+1$ and $a_{C',j}=-2a-1$, so  $\frac{-\limit-1}{2} \le a \le \frac{\limit-1}{2}$.}
\end{align*}

\paragraph{Parameters in the Blowing-up}
In three dimensions we have to consider the following two different cases of blowing-ups:
\begin{itemize}
\item When blowing-up in a full-dimensional cone we can do the same as in the two-dimensional case.\\
Denote by $C_1=\{r_1,n_1,n_2\},\ C_2=\{r_2,n_1,n_2\}$ two full-dimensional cones containing a codimension $1$ cone $C=\{n_1,n_2\}$ with $r_1+r_2 = a_{C,1}n_1+a_{C,2}n_2$. We blow up in the cone $C_2$, so the new full-dimensional cones containing $C$ are $C_1$ and $C_2'=\{r_2',n_1,n_2\}$ with $r_2' = r_2+n_1+n_2$. We obtain
\[
r_1+r_2' = a_{C,1}n_1+a_{C,2}n_2+n_1+n_2 = (a_{C,1}+1)n_1+(a_{C,2}+1)n_2
\]
and this gives $a_{C,i}+1\ge-\limit$ and thus $a_{C,i}\ge-\limit-1$. That means we have to loosen the bounds on a parameter by one if it appears at one of the codimension $1$ cones contained in the cone of the current blowing-up, just as in the two-dimensional case.
\item For the blowing-up in a codimension $1$ cone choose $C_1,\ C_2$ and $C$ as in the full-dimensional case.
When blowing up in the cone $C$ at $n=n_1+n_2$ we get six new cones $C_1'=\{r_1,n_1,n\},\ C_2'=\{r_2,n_1,n\}$ containing $C'=\{n_1,n\}$ and $C_1''=\{r_1,n_2,n\},\ C_2''=\{r_2,n_2,n\}$ containing $C''=\{n,n_2\}$. Now check the parameters at the new cone $C'$:
\[
r_1+r_2 = a_{C,1}n_1+a_{C,2}n_2 = a_{C,1}n_1 + a_{C,2}(n-n_1) = (a_{C,1}-a_{C,2})n_1 +  a_{C,2}n
\]
and at $C''$ we get equivalently:
\[
r_1+r_2 = a_{C,1}n_1+a_{C,2}n_2 = a_{C,1}(n-n_2) + a_{C,2}n_2 = a_{C,1}n +  (a_{C,2}-a_{C,1})n_2
\]
So, in the blowing-up, we have found two codimension $1$ cones each having one of the edge-parameters of the original fan, $a_{C',2}=a_{C,2}$ and $a_{C'',1}=a_{C,1}$. These edge-parameters of the two new cones $C',\ C''$ yield the same bounds for the parameters of the fan as the codimension $1$ cone $C$ in the original fan. Which means we can just pass on the bounds to the new blowing-up.
\end{itemize}

So we have again all neccessary bounds for the parameters of all fans and can create a list of completely defined smooth fans that serve as normal fans for our polytopes.

\subsubsection{Creating the Polytopes}
We can now continue as in the general and the two-dimensional case using \polymake in combination with \fourtitwo, \normaliz and \lattE to compute all smooth polytopes with at most $12$ lattice points.

\subsubsection{Results}
Running this algorithm for the remaining five minimal smooth fans yields the following result.
\begin{prop}
There are 33 smooth three-dimensional lattice polytopes with at most \nlimit lattice points.
\begin{center}
\begin{tabular}{r|cccc}
\toprule
Vertices 		& 4 & 6 & 8 & $\ge$10 \\
\hline
Polytopes 			& 2 & 25 & 6 & 0  \\
\bottomrule
\end{tabular}
\end{center}
Again the list of polytopes can be found in the appendix.
\end{prop}

\section{Final Comments}
We now have a list of smooth lattice polytopes in dimension two and three with at most $12$ lattice points. For dimension three the bound $12$ may seem rather low and in fact none of the $33$ lattice polytopes from this classification has an interior lattice point. So these polytopes may not bring great insight, but this is a first step and there are several points in the algorithm where it could be improved.\par
By implementing a different way to directly generate all smooth normal fans one could skip the big recursion calculating all blowing-ups, which is one of the most time-con\-su\-ming steps in the classification algorithm. And one could then overcome the limits of at most $12$ vertices and thus lattice points and dimension at most three.\par
The second point where one could work on is the calculation of lattice points of the polytope containing all right-hand sides. The dimension of this polytope is equal to the number of rays of the fan minus the dimension. So when we want to reach more complex smooth polytopes we will have rather high-dimensional polytopes containing the right-hand sides, so we might need better bounds for the right-hand sides and/or faster algorithms for computing the lattice points.
\pagebreak[4]
\bibliographystyle{amsplain}
\bibliography{dipl}
\appendix
{\tiny
\section{List of Smooth Polygons}
\input{polygons.tex}
\section{List of Smooth Polytopes}
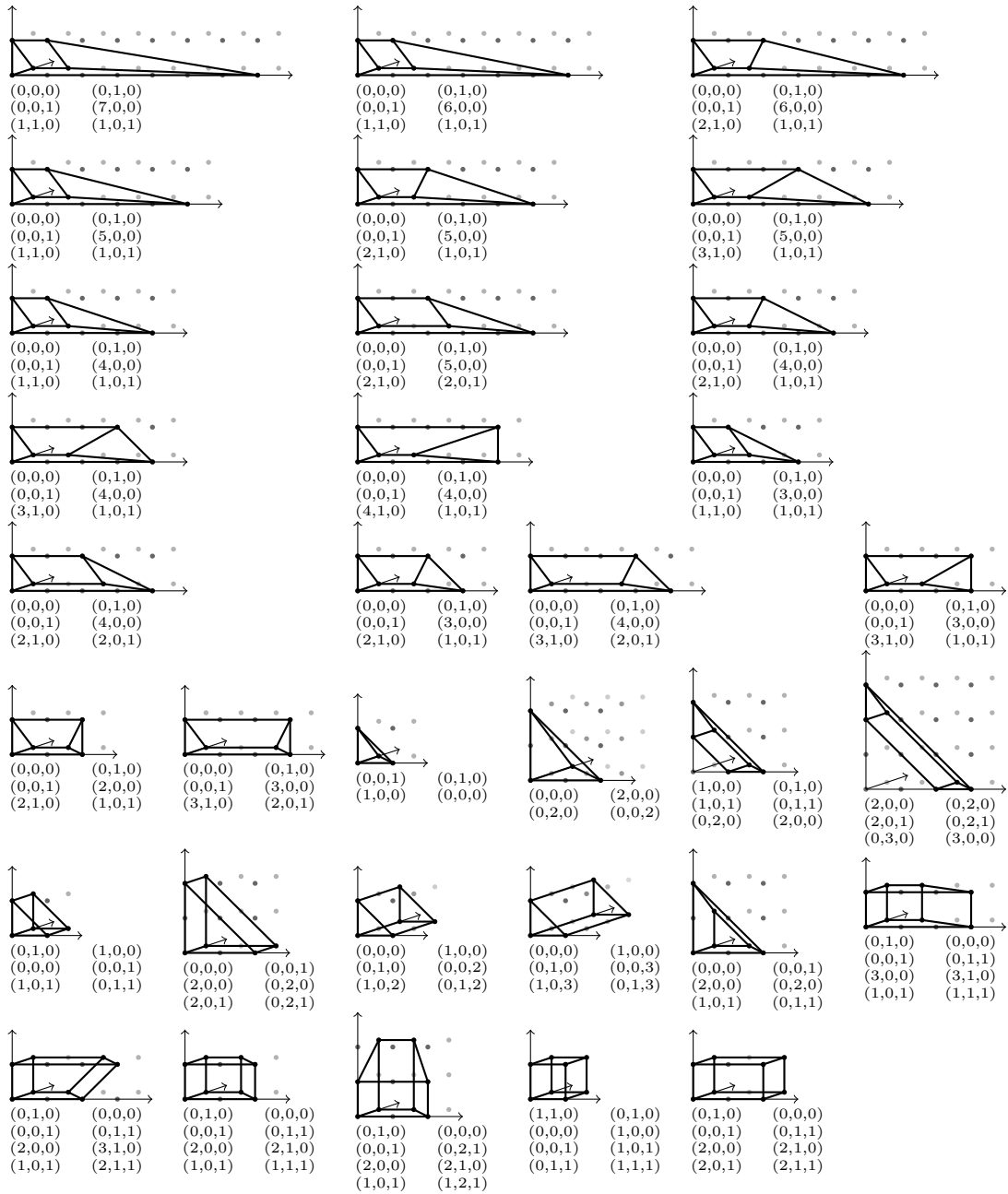
\input{polytopes.tex}
}

\pagebreak[4]

\section{\polymake extension}
\subsection{Parameterized Fans ({\ttfamily pfan.rules})}
\begin{lstlisting}[language=Perl]

# a parameterized class of fans
declare object ParameterFan<Rational> {

# maximal number of lattice points
property POINT_LIMIT : Int;

# dimension of the fan
property DIM : Int;

# incidence sets (of the rays) defining the full-dimensional cones or 
# vertices of the polytope
property VERTICES : Array<Set>;

# number of vertices
property N_VERTICES : Int;

# incidence sets (of the vertices) defining the codimension 1 cones or 
# edges between vertices of the polytopes
property EDGES : Array<Set>;

# number of edges
property N_EDGES : Int;

# rays of the fan
property RAYS : Matrix<Rational>;

# number of rays
property N_RAYS : Int;

# number of different parameters
property N_PARAMETER : Int;

#  rows with positions in RAYS matrix, parameter number and coefficient
property PARAMETER_POS : Matrix<Int>;

# lower and upper bound for each parameter (rowwise)
property PARAMETER_BOUNDS : Matrix<Int>;

# searches full-dimensional cones sharing a common codimension 1 cone
rule EDGES : VERTICES , N_VERTICES , DIM {
  my $vert = $this->VERTICES;
  my @edges = ();
  local $_;
  # check all possible pairs
  foreach (all_subsets_of_k(2,0..($this->N_VERTICES-1))) {
    my $r = $vert->[$_->[0]] * $vert->[$_->[1]];
    if ($#$r == ($this->DIM - 2)) {
      push @edges, new Set($_);
    }
  }
  $this->EDGES = new Array<Set>(@edges);
}

rule DIM : RAYS {
  $this->DIM = $this->RAYS->cols;
}

rule N_EDGES : EDGES {
  $this->N_EDGES = $#{$this->EDGES} + 1;
}

rule N_VERTICES : VERTICES {
  $this->N_VERTICES = $#{$this->VERTICES} + 1;
}

rule N_RAYS : RAYS {
  $this->N_RAYS = $this->RAYS->rows;
}

# main method  creating the list of polytopes for this fan
# calls all necessary other functions 
# returns:       the polytopes as perl array of Polytope objects
# input:           integer value determining the print- or debug-level
user_method POLYTOPES {
  my ($this,$print) = @_;
  my @poly = ();
  print "Creating blowing-ups ...\n" if($print);
  my @blowups = $this->create_blowups($print);
  print "Done creating blowing-ups.\n" if($print);
  my $i=1;
  print "Generating polytopes ...\n" if($print);
  foreach my $pfan (@blowups) {
    print "blowup ",$i," of ",scalar(@blowups),"\n" if ($print);
    my @fans = $pfan->getfans();
    my $j = 1;
    foreach my $fan (@fans) {
      print "\tfan ",$i,":",$j," of ",scalar(@blowups),":",
                scalar(@fans),"\n" if ($print);
      my $num = $fan->N_POLYTOPES;
      if ($num > 0) {
        push @poly, $fan->POLYTOPE_LIST;
      }
      $j++;
    }
    $i++;
  }
  print "Done generating polytopes.\n" if($print);
  print "Found ",scalar(@poly), " polytopes including duplicates.\n" if($print);
  print "Removing duplicates ...\n" if($print);
  my @final_poly = ();
  outer: for (my $i = scalar(@poly)-1; $i >= 0; $i--) {
    my $p = @poly[$i];
    for (my $j = $i-1; $j >= 0; $j--) {
      my $q = @poly[$j];
      # reduce neccessary isomorphy checks by checking number of vertices first
      if ($p->N_VERTICES == $q->N_VERTICES) {
        if (lattice_isomorphic_smooth_polytopes($p,$q)) {
          next outer;
        }
      }
    }
    push @final_poly, $p;
  }

  print "Found ",scalar(@final_poly), " non-isomorphic polytopes with this normal fan.\n";
  return @final_poly;
}

# starts the recursion generating all blowing-ups returning them in an perl array
# printlimit specifies up to which depth to print messages
user_method create_blowups {
  my ($fan,$printlimit) = @_;
  # incidence vectors storing which vertices and edges should be blown-up in the tree, 
  # zero vectors at the beginning
  my $vertexlist = new Vector<Int>($fan->POINT_LIMIT);
  my $edgelist = undef;
  $edgelist = new Vector<Int>(($fan->DIM * $fan->POINT_LIMIT)/2) if ($fan->DIM > 2);
  
  my @blowups = recursion($fan,$vertexlist,$edgelist,0,$printlimit);
  return @blowups;
}

# creates a new fan for every combination of parameter values and 
# returns the list of fans as perl array
user_method getfans {
  my $bf = shift;
  my $rays = $bf->RAYS;
  my @fans = ();
  if ($bf->N_PARAMETER > 0) {
    my $npar = $bf->N_PARAMETER;
    my $pb = $bf->PARAMETER_BOUNDS;
    # vector storing the current parameter values
    my $par = new Vector<Int>($pb->col(0));
    my $cont = 1;
    do {
      my $r = new Matrix<Rational>($rays);
      # add the parameter values to every position it occurs
      foreach my $pos (@{$bf->PARAMETER_POS}) {
        $r->($pos->[0],$pos->[1]) += $par->[$pos->[2]] * $pos->[3];
      }
      push @fans,new Fan<Rational>(POINT_LIMIT=>($bf->POINT_LIMIT),DIM=>($bf->DIM),
            RAYS=>$r,VERTICES=>($bf->VERTICES),EDGES=>($bf->EDGES));
      
      # generate next vector of parameter values
      # increase first parameter and check bounds,...
      $par->[0] ++;
      for (my $i = 0; $i < $npar; $i++) {
        if($par->[$i] > $pb->[$i]->[1]) {
          if($i == $npar-1) {
            $cont = 0;
          } else {
            $par->[$i] = $pb->[$i]->[0];
            $par->[$i+1]++;
          }
        }
      }
    } while($cont == 1);
  } else {
    push @fans,new Fan<Rational>(POINT_LIMIT=>($bf->POINT_LIMIT),DIM=>($bf->DIM),
          RAYS=>$rays,VERTICES=>($bf->VERTICES),EDGES=>($bf->EDGES));
  }
  return @fans;
}

# main recursion creating the blowing-ups
sub recursion {
  my ($fan,$vertexlist,$edgelist,$depth,$printlimit) = @_;
  my @b = ();
  
  if ($fan->DIM > 2) {
    # in dimension > 2 check smooth polygon criterion
    if (checkpolygons($fan) <= $fan->POINT_LIMIT) {
      push @b, $fan;
    }
  } else {
    push @b, $fan;
  }
  
  # if blowing-up does not exceed maximal number of vertices
  if ($fan->N_VERTICES <= $fan->POINT_LIMIT - ($fan->DIM-1)) {
    # blowing-up of full-dimensional cones
    for(my $i = 0; $i < $fan->N_VERTICES; $i++) {
      if ($vertexlist->[$i] < 1) {
        print "do vertexblowup:",$depth,"-",$i,"\n"if ($depth < $printlimit);
        my ($newfan,$elist) = do_vertex_blowup($fan,$i,$edgelist);
        push @b, recursion($newfan,(new Vector<Int>($vertexlist)),
                                $elist,$depth+1,$printlimit);
        $vertexlist->[$i]=1;
      }
    }
    # if dim > 2 blowing-up of codim 1 cones
    if($fan->DIM > 2) {
      for (my $i = 0; $i < $fan->N_EDGES; $i++) {
        if ($edgelist->[$i] < 1) {
          print "do edgeblowup:",$depth,"-",$i,"\n" if ($depth < $printlimit);
          my ($newfan,$vlist,$elist) = do_edge_blowup($fan,$i,$vertexlist,$edgelist);
          push @b, recursion($newfan,$vlist,$elist,$depth+1,$printlimit);
          if($edgelist->[$i] == 0) {
            $edgelist->[$i] = 1;
          }
        }
      }
    }
  }
  return @b;
}

# create the new fan after the blowing-up of full-dim cone $i
sub do_vertex_blowup {
  my ($fan,$i,$edgelist) = @_;
  # set of rays defining the vertex
  my $list = $fan->VERTICES->[$i];
  my $nrays = $fan->N_RAYS;
  my $elist = undef;
  $elist = new Vector<Int>($edgelist) if (defined($edgelist));
  
  # create and add the new ray
  my $ray = new Vector<Rational>($fan->DIM);
  local $_;
  foreach (@$list) {
    $ray += $fan->RAYS->[$_];
  }
  my $newrays = $fan->RAYS / $ray;
  
  # prepare edges array, vertices array and array with positions for new vertices
  my @newedges = @{$fan->EDGES};
  my @pos = ($i);
  my @vertices = @{$fan->VERTICES};
  for (my $k = 0; $k < $fan->DIM - 1; $k++) {
    push @pos, ($fan->N_VERTICES+$k);
    push @vertices, new Set();
  }
  local $_;
  foreach (all_subsets_of_k(2,@pos)) {
    push @newedges, (new Set($_));
  }
  
  my $k = 0;
  my $newvertices = new Array<Set>(@vertices);
  my $newedges = new Array<Set>(@newedges);
  # find neighbors of vertex we are cutting off, correct the edges 
  # and create new vertices
  for (my $j = 0; $j < $fan->N_EDGES; $j++) {
    my $edge = $fan->EDGES->[$j];
    if ($edge->contains($i))  {
      my $v = ($edge - $i)->[0];
      # new vertex
      $newvertices->[@pos[$k]] = ($fan->VERTICES->[$v] * $list) + $nrays;
      # set new edge endpoints
      $newedges->[$j] = new Set($v,@pos[$k]);
      # set neighboring edges to always ignore, 
      # we do vertex -> neighboring edge via edge -> edge
      $elist->[$j] = 2 if (defined($edgelist));
      $k++;
    }
  }
  
  my $newfan;
  # add parameters for the new ray and loosen bounds if neccessary 
  my $par = new Matrix<Int>();
  if ($fan->N_PARAMETER > 0) {
    my @param = @{$fan->PARAMETER_POS};
    local $_;
    my $p_indices = new Vector<Int>($fan->N_PARAMETER);
    # find parameters appearing in the rays defining the currently blown-up cone
    foreach (@{$fan->PARAMETER_POS}) {
      if ($list->contains($_->[0])) {
        push @param, (new Vector<Int>($nrays,$_->[1],$_->[2],$_->[3]));
        $p_indices->[$_->[2]]=1;
      }
    }
    $par = new Matrix<Int>(@param);
    # loosen bounds on all parameters found earlier
    my $pbounds = new Matrix<Int>($fan->PARAMETER_BOUNDS);
    for (my $j = 0; $j < $fan->N_PARAMETER; $j++) {
      if ($p_indices->[$j] > 0) {
        $pbounds->[$j]->[0] --;
        $pbounds->[$j]->[1] ++;
      }
    }
    # create new fan with all modified data
    $newfan = new ParameterFan<Rational>(POINT_LIMIT=>$fan->POINT_LIMIT,
        DIM=>$fan->DIM,VERTICES=>$newvertices,EDGES=>$newedges,RAYS=>$newrays,
        PARAMETER_POS=>$par,PARAMETER_BOUNDS=>$pbounds,N_PARAMETER=>$fan->N_PARAMETER);
  } else {
    $newfan = new ParameterFan<Rational>(POINT_LIMIT=>$fan->POINT_LIMIT,
        DIM=>$fan->DIM,VERTICES=>$newvertices,EDGES=>$newedges,RAYS=>$newrays,
        PARAMETER_POS=>$par,N_PARAMETER=>$fan->N_PARAMETER);
  }
  return ($newfan,$elist);
}

# create the new fan after the blowing-up of full-dim cone $i
sub do_edge_blowup {
  my ($fan,$i,$vertexlist,$edgelist) = @_;
  my $vlist = new Vector<Int>($vertexlist);
  my $elist = new Vector<Int>($edgelist);
  # indices of the two vertices of the edge
  my $edge = $fan->EDGES->[$i];
  my $e1 = $edge->[0];
  my $e2 = $edge->[1];
  # sets of the rays defining these vertices
  my $v1 = $fan->VERTICES->[$e1];
  my $v2 = $fan->VERTICES->[$e2];
  # two rays defining this edge (as set and separate)
  my $er = $v1*$v2;
  my $er1 = $er->[0];
  my $er2 = $er->[1];
  # total number of vertices and rays
  my $nv = $fan->N_VERTICES;
  my $nrays = $fan->N_RAYS;

  # add new ray
  my $newrays = $fan->RAYS / ($fan->RAYS->[$er->[0]] + $fan->RAYS->[$er->[1]]);
  
  # new edgearray
  my @newedges = @{$fan->EDGES};
  # create edges to and between the two vertices that will be put at the end
  push @newedges, new Set($nv,$nv+1);
  push @newedges, new Set($e1,$nv);
  push @newedges, new Set($e2,$nv+1);
  my $k = 0;
  
  # create new vertices
  my $newvertices = new Array<Set>(@{$fan->VERTICES},new Set(), new Set());
  $newvertices->[$e1] = $v1 - $er2 + $nrays;
  $newvertices->[$e2] = $v2 - $er2 + $nrays;
  $newvertices->[$nv] = $v1 - $er1 + $nrays;
  $newvertices->[$nv+1] = $v2 - $er1 + $nrays;
  
  # set the new vertices  to 'consider' for later blowing-ups
  $vlist->[$e1] = 0;
  $vlist->[$e2] = 0;
  
  # create the new edges and update all old edges
  # we skip all just added edges and the cut-off edge
  my $newedges = new Array<Set>(@newedges);
  for (my $j = 0; $j < $fan->N_EDGES; $j++) {
    next if ($i==$j);
    my $e = $fan->EDGES->[$j];
    if ($e->contains($e1) || $e->contains($e2)) {
      my $v = ($e - $edge)->[0];
      if (!($fan->VERTICES->[$v]->contains($er1))) {
        $newedges->[$j] = new Set($v,$nv) if ($e->contains($e1));
        $newedges->[$j] = new Set($v,$nv+1) if ($e->contains($e2));
      }
      $elist->[$j] = -1 if ($elist->[$j]==1);
    }
  }
  
  my $newfan;
  # check if the new ray has parameters
  my $par = new Matrix<Int>();
  if ($fan->N_PARAMETER > 0) {
    my @param = @{$fan->PARAMETER_POS};
    local $_;
    # search for parameters for the rays defining the new one
    foreach (@{$fan->PARAMETER_POS}) {
      if ($er->contains($_->[0])) {
        push @param, (new Vector<Int>($nrays,$_->[1],$_->[2],$_->[3]));
      }
    }
    $par = new Matrix<Int>(@param);
    $newfan = new ParameterFan<Rational>(POINT_LIMIT=>$fan->POINT_LIMIT,
            DIM=>$fan->DIM,VERTICES=>$newvertices,EDGES=>$newedges, RAYS=>$newrays,
            PARAMETER_POS=>$par,PARAMETER_BOUNDS=>$fan->PARAMETER_BOUNDS,
            N_PARAMETER=>$fan->N_PARAMETER);
  } else {
    $newfan = new ParameterFan<Rational>(POINT_LIMIT=>$fan->POINT_LIMIT,
            DIM=>$fan->DIM,VERTICES=>$newvertices,EDGES=>$newedges, RAYS=>$newrays,
            N_PARAMETER=>$fan->N_PARAMETER);
  }
  return ($newfan,$vlist,$elist);
}

# check the smooth polygon lower bound for the number of lattice points
sub checkpolygons{
  my ($fan) = @_;
  my $rays = $fan->N_RAYS;
  # results from the 2-dimensional classification, 
  # "13" means not possible for <= 12 lattice points
  # FIXME for more than POINT_LIMIT > 12 
  my $innercounts  = new Vector<Rational>('0 0 1 1 13 4 13 13');
  my $bordercounts = new Vector<Rational>('3 4 7 6 13 8 13 13');
  
  # collect degrees of all rays (number of vertices they appear in)
  my $sizes = new Vector<Int>( '0 ' x $rays );
  local $_;
  foreach (@{$fan->VERTICES}) {
    foreach my $x (@$_) {
      $sizes->[$x] ++;
    }
  }
  
  # find maximal degree, and sum up all interior lattice points
  my $count = 0;
  my $max = 0;
  local $_;
  foreach (@$sizes) {
    if ($_ > $max) {
      $max = $_;
    }
    next if ($_ > 8);
    $count += $innercounts->[$_-3];
  }
  if ($max > 8) {
    return 13;
  }
  # add boundary lattice points for the largest polygon
  $count += $bordercounts->[$max-3];
  $count += $fan->N_VERTICES - $max;
  return $count;
}

}
\end{lstlisting}
\subsection{Fans ({\ttfamily fan.rules})}
\begin{lstlisting}[language=Perl]

# a completely defined fan
declare object Fan<Rational> {

# maximal number of lattice points
property POINT_LIMIT : Rational;

# dimension of the fan
property DIM : Int;

# incidence sets (of the rays) defining the full-dimensional cones or 
# vertices of the polytope
property VERTICES : Array<Set>;

# number of vertices
property N_VERTICES : Int;

# incidence sets (of the vertices) defining the codimension 1 cones or 
# edges between vertices of the polytopes
property EDGES : Array<Set>;

# number of edges
property N_EDGES : Int;

# rays of the fan
property RAYS : Matrix<Rational>;

# number of rays
property N_RAYS : Int;

# polytope approximating the set of valid right-hand sides
property RHS_POLYTOPE : Polytope<Rational>;

# set of right-hand side vectors defining polytopes 
# with at most POINT_LIMIT lattice points
property RHS_VECTORS : Set<Vector<Integer>>;

# set of right-hand side vectors defining polytopes 
# with at most POINT_LIMIT lattice points
property N_POLYTOPES : Integer;

# returns a perl array of all polytopes with this normal fan and 
# at most POINT_LIMIT lattice points
user_method POLYTOPE_LIST {
  my $fan = shift;
  my @list = ();
  # create polytope for every rhs vector
  foreach my $rhs (@{$fan->RHS_VECTORS}) {
    my $ineq = (new Vector<Rational>($rhs)) | ($fan->RAYS);
    my $p = new Polytope<Rational>(INEQUALITIES=>$ineq);
    #my $v = $p->N_VERTICES;
    #my $lp = $p->N_LATTICE_POINTS;
    push @list, $p;
  }
  return @list;
}

# searches full-dimensional cones sharing a common codimension 1 cone
rule EDGES : VERTICES , N_VERTICES , DIM {
  my $vert = $this->VERTICES;
  my @edges = ();
  local $_;
  # check all possible pairs
  foreach (all_subsets_of_k(2,0..($this->N_VERTICES-1))) {
    my $r = $vert->[$_->[0]] * $vert->[$_->[1]];
    if ($#$r == ($this->DIM - 2)) {
      push @edges, new Set($_);
    }
  }
  $this->EDGES = new Array<Set>(@edges);
}

rule DIM : RAYS {
  $this->DIM = $this->RAYS->cols;
}

rule N_EDGES : EDGES {
  $this->N_EDGES = $#{$this->EDGES} + 1;
}

rule N_VERTICES : VERTICES {
  $this->N_VERTICES = $#{$this->VERTICES} + 1;
}

rule N_RAYS : RAYS {
  $this->N_RAYS = $this->RAYS->rows;
}

rule N_POLYTOPES : RHS_VECTORS {
  $this->N_POLYTOPES = $#{$this->RHS_VECTORS}+1;
}

# searches in RHS_POLYTOPE for RHS-vectors defining polytopes 
# with at most POINT_LIMIT lattice points
rule RHS_VECTORS : RHS_POLYTOPE , RAYS , POINT_LIMIT , DIM {
  my $valid = $this->RHS_POLYTOPE;
  my $limit = $this->POINT_LIMIT;
  my $dim = $this->DIM;
  my $rays = $this->RAYS;
  
  # check the polytope for feasibility and number of lattice points first
  if ($valid->FEASIBLE == 0) {
    $this->RHS_VECTORS=new Set<Vector<Integer>>();
    return;
  }
  if ($valid->N_LATTICE_POINTS == 0) {
    $this->RHS_VECTORS=new Set<Vector<Integer>>();
    return;
  }
  
  # remove first coordinate (homogenized)
  my $rhs = $valid->LATTICE_POINTS->minor(All,range(1,$rays->rows()-$dim));
  my @rhs = @$rhs;
  @rhs = sort @rhs;
  
  my @good_rhs = ();
  my $i=0;
  my $sched2;
  
  while (@rhs > 0) {
    # create new RHS vector and inequality matrix
    my $x = ((zero_vector<Int>($dim)) | shift(@rhs));
    my $ineq = (new Vector<Rational>($x)) | $rays;
    
    # define polytope and check number of lattice points
    my $p = new Polytope<Rational>(INEQUALITIES=>$ineq);
    
    if ($p->N_LATTICE_POINTS <= $limit) {
      # found a new polytope
      push @good_rhs, $x;
      $i++;
    } else {
      # we only keep RHS-vectors with at least one entry with smaller 
      # absolute value than the current vector
      my @newrhs = ();
      point: for my $v (@rhs) {
        for(my $j = 0; $j < $rays->rows()-$dim; $j++) {
          if (abs($v->[$j]) < abs($x->[$j+$dim])) {
            push @newrhs,$v;
            next point;
          }
        }
      }
      @rhs = @newrhs;
    }
  }
  $this->RHS_VECTORS=new Set<Vector<Integer>>(@good_rhs);
}

# creates the polytope approximating the set of right-hand sides
rule RHS_POLYTOPE : RAYS , VERTICES , EDGES , POINT_LIMIT , N_EDGES , N_VERTICES {
  my $vert = $this->VERTICES;
  my $edges = $this->EDGES;
  my $rays = $this->RAYS;
  my $dim = $rays->cols();
  my $r = $rays->rows();
  my $limit = $this->POINT_LIMIT;
  
  # stores the sum of all edge-lengths
  my $totallength = new Vector<Rational>($r+1);
  
  # first stores l(e) >= 1 inequalities, second one inequalities for l(e)<=...
  my @ineq = ();
  my @l_ineq = ();
  
  # calculate all edge lengths and create inequalities
  foreach (@$edges) {
    # vertices defining the edge
    my $fromvert = $vert->[$_->[0]];
    my $tovert = $vert->[$_->[1]];
    
    # rays defining the edge (n_1,...,n_{d-1})
    my $neighbors = $fromvert * $tovert;
    
    # rays limiting this edge (r_1,r_2)
    my $from = $fromvert - $neighbors;
    my $to = $tovert - $neighbors;
    $from = $from->[0];
    $to = $to->[0];
    
    # calculate inequalities
    my ($edgelimit,$pointlimit) = edgelimits($rays,$from,$to,$limit,$neighbors);
    
    $totallength += $pointlimit;
    push @ineq, $edgelimit;
    push @l_ineq, $pointlimit;
  }
  
  # right-hand side for total length inequality
  $totallength->[0] = $limit + $this->N_EDGES - $this->N_VERTICES;
  push @l_ineq, $totallength;
  
  my $ineq = new Matrix<Rational>(\@ineq);
  my $l_ineq = new Matrix<Rational>(\@l_ineq);
  
  # for all rays with all entries positive we know that the right-hand side must 
  # be postive and vice versa (we want the negative RHS-vectors)
  my $signs = new Matrix<Rational>($r,$r+1);
  my @signs = ();
  row: for (my $i = 0; $i < $r; $i++) {
    my $sign = 0;
    for (my $j = 0; $j < $dim; $j++) {
      if ($rays->($i,$j) < 0) {
        next row if ($sign > 0);
        $sign = -1;
      } elsif ($rays->($i,$j) > 0) {
        next row if ($sign < 0);
        $sign = 1;
      }
    }
    $signs->($i,$i+1) = -$sign;
  }
  
  # put all inequalities together and remove the first d columns after the rhs entry
  # because those variables are always zero (we do not want translations)
  $ineq = $ineq / $signs;
  $ineq = $ineq->minor(All,range(0,0)+range($dim+1,$r));
  $l_ineq = $l_ineq->minor(All,range(0,0)+range($dim+1,$r));
  my $limited = (new Matrix<Rational>($ineq)) / (new Matrix<Rational>($l_ineq));  
  
  $this->RHS_POLYTOPE = new Polytope<Rational>( INEQUALITIES=>$limited );
}

# calculates the inequalities in terms of the rhs vector given by 1<=l(e)<=...
sub edgelimits {
  my ($rays, $from, $to, $limit, $nb) = @_;
  my $dim = $rays->cols();
  # limiting rays
  my $r1 = $rays->row($to);
  my $r2 = $rays->row($from);
  # rays defining the edge
  my $N = $rays->minor($nb,All);
  
  # determine edge-parameters a_{C,i}
  $N = new Matrix<Rational>(transpose($N));
  my $v = new Vector<Rational>($r1+$r2);
  my $a = lin_solve($N,$v);
  my $sum = 0;
  # and sum them up for the right-hand side
  $sum += $a->[$_] for (0 .. ($dim-2));
  
  # vectors storing the edgelength >= 1 and <=... inequalities
  my $edgelimit = new Vector<Rational>($rays->rows() + 1);
  my $pointlimit = new Vector<Rational>($rays->rows() + 1);
  
  # rhs of the inequalities
  $edgelimit->[0] = (-1);
  $pointlimit->[0] = (($limit - $sum)/$dim-1);
  
  # set coefficients
  for (my $i = 0; $i < $dim-1; $i++) {
    my $k = $nb->[$i] + 1;
    $edgelimit->[$k] = (-1) * $a->[$i];
    $pointlimit->[$k] = $a->[$i];
  }
  $edgelimit->[$to+1] = 1;
  $pointlimit->[$to+1] = -1;
  $edgelimit->[$from+1] = 1;
  $pointlimit->[$from+1] = -1;
  
  return ($edgelimit,$pointlimit);
}

}
\end{lstlisting}
\subsection{Additional Script ({\ttfamily checkiso})}
\begin{lstlisting}[language=Perl]
# removes all lattice isomorphic polytopes from a list
# only the first one of an isomorphy class is kept
sub checkiso {
  my (@list) = @_;
  my @res = ();
  outer: for (my $i = scalar(@list)-1; $i >= 0; $i--) {
    my $p = @list[$i];
    for (my $j = $i-1; $j >= 0; $j--) {
      my $q = @list[$j];
      if ($p->N_VERTICES == $q->N_VERTICES) {
        if ($p->N_LATTICE_POINTS == $q->N_LATTICE_POINTS) {
          next outer if (lattice_isomorphic_smooth_polytopes($p,$q));
        }
      }
    }
    push @res, $p;
  }
  return @res;
}
\end{lstlisting}

\subsection{Minimal Fans as \polymake files}
\subsubsection{Dimension Two}
\paragraph{Projective Space}
\begin{lstlisting}[language=xml]
<?xml version="1.0" encoding="utf-8"?>

<object name="P" type="polytope::ParameterFan&lt;Rational&gt;" version="2.9.7" 
tm="4b0d8c7b" xmlns="http://www.math.tu-berlin.de/polymake/#3">
  <property name="POINT_LIMIT" value="12" />
  <property name="DIM" value="2" />
  <property name="VERTICES">
    <m>
      <v>0 1</v>
      <v>1 2</v>
      <v>0 2</v>
    </m>
  </property>
  <property name="RAYS">
    <m>
      <v>1 0</v>
      <v>0 1</v>
      <v>-1 -1</v>
    </m>
  </property>
  <property name="N_RAYS" value="3" />
  <property name="N_VERTICES" value="3" />
  <property name="N_PARAMETER" value="0" />
</object>
\end{lstlisting}
\paragraph{Hirzebruch surfaces}
\begin{lstlisting}[language=xml]
<?xml version="1.0" encoding="utf-8"?>

<object name="H" type="polytope::ParameterFan&lt;Rational&gt;" version="2.9.7" 
tm="4af96f7b" xmlns="http://www.math.tu-berlin.de/polymake/#3">
  <property name="POINT_LIMIT" value="12" />
  <property name="DIM" value="2" />
  <property name="VERTICES">
    <m>
      <v>0 1</v>
      <v>0 2</v>
      <v>1 3</v>
      <v>2 3</v>
    </m>
  </property>
  <property name="RAYS">
    <m>
      <v>1 0</v>
      <v>0 1</v>
      <v>0 -1</v>
      <v>-1 0</v>
    </m>
  </property>
  <property name="N_RAYS" value="4" />
  <property name="N_VERTICES" value="4" />
  <property name="PARAMETER_POS">
    <m>
      <v>3 1 0 1</v>
    </m>
  </property>
  <property name="PARAMETER_BOUNDS">
    <m>
      <v>0 12</v>
    </m>
  </property>
  <property name="N_PARAMETER" value="1" />
</object>
\end{lstlisting}
\subsubsection{Dimension Three}
\paragraph{Fan $\mathbf{3^4}$}
\begin{lstlisting}[language=xml]
<?xml version="1.0" encoding="utf-8"?>

<object name="f" type="polytope::ParameterFan&lt;Rational&gt;" version="2.9.7" 
tm="4b0d8c1d" xmlns="http://www.math.tu-berlin.de/polymake/#3">
  <property name="POINT_LIMIT" value="12" />
  <property name="DIM" value="3" />
  <property name="VERTICES">
    <m>
      <v>0 1 2</v>
      <v>0 1 3</v>
      <v>0 2 3</v>
      <v>1 2 3</v>
    </m>
  </property>
  <property name="RAYS">
    <m>
      <v>1 0 0</v>
      <v>0 1 0</v>
      <v>0 0 1</v>
      <v>-1 -1 -1</v>
    </m>
  </property>
  <property name="N_RAYS" value="4" />
  <property name="N_VERTICES" value="4" />
  <property name="N_PARAMETER" value="0" />
</object>
\end{lstlisting}
\paragraph{Fan $\mathbf{(3^24^2)'}$}
\begin{lstlisting}[language=xml]
<?xml version="1.0" encoding="utf-8"?>

<object name="f" type="polytope::ParameterFan&lt;Rational&gt;" version="2.9.7" 
tm="4b0d571c" xmlns="http://www.math.tu-berlin.de/polymake/#3">
  <property name="POINT_LIMIT" value="12" />
  <property name="DIM" value="3" />
  <property name="VERTICES">
    <m>
      <v>0 1 2</v>
      <v>0 1 3</v>
      <v>0 2 4</v>
      <v>0 3 4</v>
      <v>1 2 4</v>
      <v>1 3 4</v>
    </m>
  </property>
  <property name="RAYS">
    <m>
      <v>1 0 0</v>
      <v>0 1 0</v>
      <v>0 0 1</v>
      <v>0 0 -1</v>
      <v>-1 -1 0</v>
    </m>
  </property>
  <property name="N_RAYS" value="5" />
  <property name="N_VERTICES" value="6" />
  <property name="PARAMETER_POS">
    <m>
      <v>4 2 0 -1</v>
    </m>
  </property>
  <property name="PARAMETER_BOUNDS">
    <m>
      <v>0 12</v>
    </m>
  </property>
  <property name="N_PARAMETER" value="1" />
</object>
\end{lstlisting}
\paragraph{Fan $\mathbf{(3^24^2)''}$}
\begin{lstlisting}[language=xml]
<?xml version="1.0" encoding="utf-8"?>

<object name="f" type="polytope::ParameterFan&lt;Rational&gt;" version="2.9.7" 
tm="4b0d8bd2" xmlns="http://www.math.tu-berlin.de/polymake/#3">
  <property name="POINT_LIMIT" value="12" />
  <property name="DIM" value="3" />
  <property name="VERTICES">
    <m>
      <v>0 1 2</v>
      <v>0 1 3</v>
      <v>0 2 3</v>
      <v>1 2 4</v>
      <v>1 3 4</v>
      <v>2 3 4</v>
    </m>
  </property>
  <property name="RAYS">
    <m>
      <v>1 0 0</v>
      <v>0 1 0</v>
      <v>0 0 1</v>
      <v>0 -1 -1</v>
      <v>-1 0 0</v>
    </m>
  </property>
  <property name="N_RAYS" value="5" />
  <property name="N_VERTICES" value="6" />
  <property name="PARAMETER_POS">
    <m>
      <v>4 1 0 1</v>
      <v>4 2 1 1</v>
    </m>
  </property>
  <property name="PARAMETER_BOUNDS">
    <m>
      <v>-12 12</v>
      <v>-12 12</v>
    </m>
  </property>
  <property name="N_PARAMETER" value="2" />
</object>
\end{lstlisting}
\paragraph{Fan $\mathbf{4^6}$}
\begin{lstlisting}[language=xml]
<?xml version="1.0" encoding="utf-8"?>

<object name="f" type="polytope::ParameterFan&lt;Rational&gt;" version="2.9.7" 
tm="4b0d8a3c" xmlns="http://www.math.tu-berlin.de/polymake/#3">
  <property name="POINT_LIMIT" value="12" />
  <property name="DIM" value="3" />
  <property name="VERTICES">
    <m>
      <v>0 1 2</v>
      <v>0 1 3</v>
      <v>0 2 4</v>
      <v>0 3 4</v>
      <v>1 2 5</v>
      <v>1 3 5</v>
      <v>2 4 5</v>
      <v>3 4 5</v>
    </m>
  </property>
  <property name="RAYS">
    <m>
      <v>1 0 0</v>
      <v>0 1 0</v>
      <v>0 0 1</v>
      <v>0 0 -1</v>
      <v>0 -1 0</v>
      <v>-1 0 0</v>
    </m>
  </property>
  <property name="N_RAYS" value="6" />
  <property name="N_VERTICES" value="8" />
  <property name="PARAMETER_POS">
    <m>
      <v>5 1 0 -1</v>
      <v>4 2 1 1</v>
      <v>5 2 2 1</v>
    </m>
  </property>
  <property name="PARAMETER_BOUNDS">
    <m>
      <v>-12 12</v>
      <v>-12 12</v>
      <v>-12 12</v>
    </m>
  </property>
  <property name="N_PARAMETER" value="3" />
</object>
\end{lstlisting}
\paragraph{Fan $\mathbf{3^24^36^2}$}
\begin{lstlisting}[language=xml]
<?xml version="1.0" encoding="utf-8"?>

<object name="f" type="polytope::ParameterFan&lt;Rational&gt;" version="2.9.7" 
tm="4b0d42cd" xmlns="http://www.math.tu-berlin.de/polymake/#3">
  <property name="POINT_LIMIT" value="12" />
  <property name="DIM" value="3" />
  <property name="VERTICES">
    <m>
      <v>0 1 2</v>
      <v>0 1 6</v>
      <v>0 2 5</v>
      <v>0 3 4</v>
      <v>0 3 5</v>
      <v>0 4 6</v>
      <v>1 2 6</v>
      <v>2 5 6</v>
      <v>3 4 6</v>
      <v>3 5 6</v>
    </m>
  </property>
  <property name="RAYS">
    <m>
      <v>1 0 0</v>
      <v>0 1 0</v>
      <v>0 0 1</v>
      <v>0 0 -1</v>
      <v>0 1 -1</v>
      <v>0 -1 0</v>
      <v>-1 2 -1</v>
    </m>
  </property>
  <property name="N_RAYS" value="7" />
  <property name="N_VERTICES" value="10" />
  <property name="PARAMETER_POS">
    <m>
      <v>5 2 0 -1</v>
    </m>
  </property>
  <property name="PARAMETER_BOUNDS">
    <m>
      <v>-6 5</v>
    </m>
  </property>
  <property name="N_PARAMETER" value="1" />
</object>
\end{lstlisting}
\end{document}